\documentclass[12pt,leqno]{article}

 \usepackage{color}

\usepackage{amsmath}
\usepackage{amsfonts,amssymb,latexsym,epsfig,amsthm, amscd, a4wide}
\usepackage{color}

\makeatletter
\@addtoreset{equation}{section}
\makeatother
    \makeatletter
    \renewcommand*{\@fnsymbol}[1]{\ensuremath{\ifcase#1 \or 1\or *\or 2\or 3\or
        \mathsection\or \mathparagraph\or \|\or **\or \dagger\dagger
        \or \ddagger\ddagger \else\@ctrerr\fi}}
    \makeatother

\newcommand{\E}        {{ {\mathbb E}}}
\renewcommand{\P}      {{ {\mathbb P}}}

\newcommand{\R}        {{\mathbb R}}

\newcommand\FF{{\mathcal F}}
\newcommand{\F}  {{\mathcal F}}

\newcommand\HH{{\mathcal H}}

\newcommand{\eps}  {{\varepsilon}}

\numberwithin{equation}{section}

\newtheorem{theorem}{{Theorem}}[section]
\newtheorem{proposition}[theorem]{{Proposition}}
\newtheorem{lemma}[theorem]{Lemma}
\newtheorem{corollary}[theorem]{Corollary}

\newtheorem{remark}[theorem]{Remark}

\newcommand{\bpf}[1][Proof]{{\noindent {\sc #1 }}}

\newcommand{\epf}{{\hfill $\square$}}
\makeindex

\begin{document}
\title{The height process of a continuous state branching process with  interaction}
\author{Zenghu Li\footnote{School of Mathematical Sciences, Beijing Normal University, Beijing 100875, China,
{lizh@bnu.edu.cn}}, \'Etienne Pardoux\footnote{corresponding author} \footnote{Aix Marseille Univ, CNRS, Centrale Marseille, I2M, Marseille, France,
{etienne.pardoux@univ-amu.fr}}, Anton Wakolbinger\footnote{Institute of Mathematics, FB 12,
Goethe-University
60629 Frankfurt,
Germany,
{wakolbinger@math.uni-frankfurt.de}}}
\date{}

\maketitle

\begin{abstract}
For a generalized continuous state branching process with non-vanishing diffusion part, finite expectation  and a directed (``left-to-right'') interaction, we construct the height process of its forest of genealogical trees. The connection between this height process and the population size process  is given by an extension of the second Ray--Knight theorem. This paper generalizes earlier work of the two last authors which was restricted to the case of continuous branching mechanisms. Our approach is different from that of Berestycki et al.~\cite{BFF}. There the diffusion part of the population process was allowed to vanish, but the class of interactions was more restricted.
\end{abstract}

{\sc Keywords} {Continuous state branching process, Population dynamics with interaction, Genealogy, Height process of a random tree}

{\sc AMS Subject Classification} {Primary 60J80; 60J25; 60H10; Secondary 92D25}
\section{Introduction}\label{Intro}
The most general continuous state branching processes (CSBP's) are solutions of SDEs of the form
\begin{equation}\label{eq}
\begin{split}
Z^x_t&=x+\gamma\int_0^t Z^x_rdr+\sqrt{2\beta}\int_0^t\int_0^{Z^x_r}W(dr,du)+
\int_0^t\int_0^{Z^x_{r-}}\int_0^1z\widetilde{M}(dr,du,dz)\\
&\qquad+\int_0^t\int_0^{Z^x_{r-}}\int_1^\infty z{M}(dr,du,dz), \quad t\ge0,
\end{split}
\end{equation}
where $W(dr,du)$ is a space--time white noise, $M(dr,du,dz)$ is a Poisson Random Measure (PRM) on $(0,+\infty)^3$
with intensity $dr\, du\, \pi(dz)$ and $\widetilde{M}(dr,du,dz)=M(dr,du,dz)-dr\, du\, \pi(dz)$. The $\sigma$--finite measure $\pi$ is assumed to be such that $(z^2\wedge1)\pi(dz)$ is a finite measure on $(0,\infty)$.

We shall  assume in this paper that 
\begin{equation}\label{hyp:pi}
\beta>0,\qquad \int_0^\infty (z^2\wedge z)\pi(dz)<\infty.
\end{equation}
The assumption $\beta > 0$ will be essential to obtain a new representation of the {\em height process} (of a genealogical forest) that underlies \eqref{eq}, see Proposition 3.13 below. This approach, using tools from stochastic analysis, will  be the basis for a representation of $H$ also in the case with interaction, see \eqref{eq:H}.  The second condition in \eqref{hyp:pi} allows us  to replace the drift coefficient $\gamma$ by $-\alpha := \gamma - \int_1^\infty z\pi(dz)$, and to write the last two integrals in equation \eqref{eq} as a single integral with respect to $\widetilde{M}$, namely
\begin{equation}\label{eq:CSBP}
\begin{split}
Z^{x}_t&=x-\alpha\int_0^t Z^x_rdr+\sqrt{2\beta}\int_0^t\int_0^{Z^{x}_r}W(dr,du)\\
&\qquad+\int_0^t\int_0^{Z^{x}_{r-}}\int_{0}^\infty z\widetilde{M}(dr,du,dz),  \quad t\ge0.
\end{split}
\end{equation}
 Moreover, we shall consider a generalized CSBP,
where the linear drift $-\alpha z$ is replaced by a nonlinear drift $f(z)$, which in general destroys the  branching property, making  $Z^x$ and $Z^{x+y}-Z^x$ dependent. Specifically, we consider the collection of  SDE's, indexed by $x \ge 0$,
\begin{equation}\label{eq:Zinteract}
\begin{split}
Z^{x}_t&=x+\int_0^t f(Z^{x}_r)dr+\sqrt{2\beta}\int_0^t\int_0^{Z^{x}_r}W(dr,du)\\
&\qquad+\int_0^t\int_0^{Z^{x}_{r-}}\int_{0}^\infty z\widetilde{M}(dr,du,dz), \quad t\ge0.
\end{split}
\end{equation}
We assume 
\begin{equation}\label{hyp:f}
f\in C^1(\R_+),\ f(0)=0,\  f'(z)\le\theta,\ \text{ for all }z\in\R,
\end{equation}
for some $\theta\in\R$.
The two assumptions \eqref{hyp:pi} and \eqref{hyp:f} will be assumed to hold throughout this paper, and will not be repeated in the statements.

It follows from Theorem 2.1 in Dawson and Li \cite{DL} that equation \eqref{eq:Zinteract} has a unique strong solution.
The introduction of the term $\int_0^t\int_0^{Z^{x}_r}W(dr,du)$ to replace the more traditional 
$\int_0^t\sqrt{Z^x_r}dB_r$ is due to \cite{DL}. Its motivation is to have a unified noise driving the equation for all initial conditions $x$. In the case of linear $f$ this provides a coupling for the CSBP's with different initial conditions. We retain that same coupling here. 

Our motivation for considering the SDE \eqref{eq:Zinteract} is to model large populations with a specific form of interaction. It is shown in Dram\'e and Pardoux \cite{DP} that an appropriately renormalized sequence of branching processes with interaction converges to the solution of \eqref{eq:Zinteract}.

In this paper we want to describe the height process $(H_s)$ of a forest of genealogical trees of the population whose total mass process $(Z_t^x)$  satisfies \eqref{eq:Zinteract}. We will always write $s$ for the ``exploration time'' and $t$ for the ``real time'', so that $H_s$ can be thought as the real time at which an individual lives that is explored at time $s$. The basic building block for the construction of $H$ is a spectrally positive L\'evy process~$X$ (see subsection \ref{sec1.2}), which due to the assumption $\beta >0$ has a Brownian component.
The equations for $H$ and the accompanying  L\'evy process $X$, then with a drift, are
\begin{equation}\label{eq:H}
\begin{split}
\beta H_s&=\int_0^sf'(L^{H_r}(r))dr+\sqrt{2\beta}B_s+\int_0^s\int_0^\infty z\tilde{N}(dr,dz)-\inf_{0\le r\le s}X_r\\&\qquad-\int_0^s\int_0^\infty\left(z+\inf_{r\le u\le s}X_u-X_r\right)^+N(dr,dz)\, ,\quad s \ge 0\, ,
\end{split}
\end{equation}
where $L^t(s)$ stands for the local time accumulated by the process $H$ at level~$t$ up to time~$s$, 
\begin{align}\label{eq:X}
X_s=\int_0^sf'(L^{H_r}(r))dr+\sqrt{2\beta} B_s+\int_0^s\int_0^\infty z\tilde{N}(dr,dz), \quad s \ge 0\, ,
\end{align}
$B$ is a standard Brownian motion, 
 $N$ is a Poisson random measure on $(0,+\infty)^2$ with mean measure $dr\,\pi(dz)$ and $\tilde{N}(dr,dz)=N(dr,dz)-dr\pi(dz)$.
We shall see that in the case $f(x)=-\alpha x$, $\alpha\ge0$, our formula for $H$ is equivalent to the formulas which appear in Duquesne and Le Gall \cite{DLG}. In the general case, we solve the SDE for $H$ with the help of Girsanov's theorem. This change of measure introduces the ``local time drift'' that appears also in \eqref{eq:X} for $X$.

We note that  \eqref{eq:Zinteract} and \eqref{eq:H} go along with a natural linear (left-to-right) ordering of the (continuum of) individuals that are alive at time $t$, and corresponds to the ordering of the exploration time $s$. This results in an  individual interaction which acts in a {\em directed} way, and is compatible with the global feedback of the population size  on the population growth that is described by the function $f$.  E.g., for $f(z)=-z^2$, the population  $Z^1$  will experience less downward drift than the population   $Z^{2}-Z^1$ ; this is the effect of the directed ``trees under attack'' dynamics that was the starting point in \mbox{Le et al.} \cite{LePW} and Pardoux and Wakolbinger \cite{PW} and was related to \eqref{eq:Zinteract} by the same authors in \cite{PWa}.  The present work thus extends previous work in case of continuous CSBPs, which started with the logistic interaction $f(z)=az-bz^2$ in  \cite{LePW} and \cite{PW}, and then described more general interactions in Ba and Pardoux \cite{BP} and in Pardoux \cite{P}.

  The connection between the height process $(H_s)$ and the population with total mass $(Z^x_t)$ will be given by an extension of the second Ray--Knight theorem, Theorem \ref{RNfinal} below, which roughly speaking says that if $L^t(s)$ denotes the local time accumulated at level $t$ by the process $H$ up to time $s$, and if $S_x=\inf\{s>0, L^0(s)>x\}$, then
$\{L^t(S_x),\, t\ge 0\}$ solves the SDE \eqref{eq:Zinteract}. In fact, since we do not know a priori whether or not the process $H$ returns to $0$ often enough such that its local time at $0$ accumulates mass $x$ (or in other words whether $Z^x$ hits zero in finite time), we will rather consider the process $H$ with an additional drift $g_a$ which modifies the dynamics of $H$ above an arbitrary level $a>0$, and insures that the process $H$ return to $0$ after any time $s>0$. The intuitive reason why this works is that, due to the  fact that $X$ has independent increments, and the properties of the Poisson random measure $N$, the pieces of trajectories of $H$ which accumulate local time at levels below $a$ interact with the past of $H$ only through the drift, which is a function of the local time accumulated at the current level, so in particular it does not depend upon the behavior of $H$ while it takes values in $(a,+\infty)$ hence it does not depend upon the additional drift $g_a$. As a result, for fixed $a$, we have the Ray--Knight interpretation only on the time interval $[0,a]$.   In order  to make sure that Girsanov's theorem is applicable, we start out by replacing  $f$ by a function $f_b$ which coincides with $f$ on the interval $[0,b]$, while $f_b$ and $f'_b$ are bounded and the latter is also uniformly continuous. The limit $b\to\infty$ leads to a family of probability measures $\P^a$, $a>0$, which admits two projective limits: one of the laws of $(H,X)$ under $\P^a$ which gives a unique weak solution of \eqref{eq:H}, \eqref{eq:X}, the other one of the laws of $\{L^t_x(S_x), 0\le t\le a, x > 0\}$ under $\P^a$, which gives the Ray-Knight representation of \eqref{eq:Zinteract}.

Berestycki, Fittipaldi and Fontbona \cite{BFF} establish an extended Ray--Knight theorem in the same situation as ours, except that, while they do not restrict themselves to the case $\beta>0$, their assumptions (in their Theorem 1.2) on the nonlinear interaction $f$ are more restrictive than our hypothesis (1.5): they assume that $f$ (with $f(0)=0$) is differentiable and concave and has a non-positive and locally Lipschitz derivative. Their equation (1.5) is the analogue to our equation 
\eqref{eq:Zinteract}, with the interaction term working ``from left to right'' like in our setting. Their approach to the underlying tree--picture is, however, quite different from ours. While they translate the competition type interaction (this follows from the non positive assumption for $f'$) into a pruning procedure on the forest of trees corresponding to the CSBP, we consider 
the process $H$ as the solution of an SDE, with a drift which is $f'$ evaluated at the local time of $H$ at time $s$ and at the level $H_s$. This is an extension of the SDE for $H$ in the case without jumps, as it appears e.g. in~\cite{BP}. The new difficulty is that each jump of $Z^x$ creates a new sub--forest of trees which must be explored. As a result, $H$ is not a Markov process. It should remember at which level (i.e. time for the process $Z^x$) a forest of trees for a certain mass of population was created, and that sub--forest should be completely explored, before the height process is allowed to go below that level.

 We shall need
to consider local times of processes which are not necessarily continuous semi--martingales. This will extend the following definition: If $Y$ is  a continuous semi--martingale,  we shall denote
by $L^a(s,Y)$, or $L^a(s)$ if there is no risk of ambiguity, the local time accumulated by the process $Y$ at level~$a$ up
to time $s$, in the sense that it satisfies
\begin{equation}\label{loctime}
  L^a(s,Y)=\lim_{\eps\to0}\frac{1}{\eps}\int_0^s{\bf1}_{[a,a+\eps]}(Y_r)dr.
  \end{equation}
It then follows from the occupation times formula that for any Borel measurable \mbox{$g:\R\to\R_+$,}
\[ \int_0^sg(Y_r)dr=\int_{-\infty}^\infty g(a)L^a_sda.\]

Our approach to the interactive case is built upon a fresh look at the height processes $H$ constructed by Duquesne and Le Gall \cite{DLG},  for general CSBP's.  We will give a new representation of $H$ which allows for an extension  to the interactive case, including the corresponding Ray-Knight representation of the solution $Z^x$ of \eqref{eq:Zinteract}.

 As a matter of fact, a large part of the present paper is concerned with the linear (CSBP) case, i.e. the case where $Z^x$ solves \eqref{eq:CSBP}. In this case, thanks to the assumption $\beta >0$, the height process obeys (see formula (1.4) in \cite{DLG})
\begin{equation}\label{1-4}
 \beta H_s=|\{\overline{X}^s_r;\ 0\le r\le s\}|,
 \end{equation}
 where  $ \overline{X}^s_r:=\inf_{r\le u\le s}X_u$ and
 $|A|$ denotes the Lebesgue measure of the set $A$. The first step of our work will consist in reinterpreting that formula, in a form which will allow the generalization to a nonlinear function $f$ (i. e. to the case of interaction).

The paper is organized as follows. Section \ref{sec:RF} is very short. It  makes precise some properties of the space--time random field $\{Z^x_t,\ t\ge0, x\ge0\}$.
Section \ref{sec:csbp} considers the case without interaction. We first  establish preliminary results that are necessary for the definition of the non-Markovian term in our representation of the height process  $H$, namely the non-compensated integral w.r.t. $N$ which appears in equation \eqref{eq:H}. We then study successively the cases $\pi=0$ (no jumps), $\pi$ finite, and finally the general case where $\pi$ satisfies \eqref{hyp:pi}, and establish the Ray--Knight theorem in a way which is tailored for the subsequent extension to the interactive case. Section \ref{sec:interact} considers the case with the interaction $f$. We introduce the SDE for $H$ which  has a drift term that depends on the local time accumulated at the current height. In order to prove the Ray-Knight representation of the solution of \eqref{eq:Zinteract} in terms of the local time of $H$, we again treat successively  the cases $\pi=0$ (no jumps), $\pi$ finite, and finally the general case where $\pi$ satisfies \eqref{hyp:pi}.

\section{The population sizes as a random field}\label{sec:RF}
The population size process $\{Z^x_t,\, t,x\ge0\}$ solving \eqref{eq:Zinteract}  is an $\R_+$--valued random field indexed by  $t$ and $x$. For each fixed $x>0$, $\{Z^x_t,\, t\ge0\}$ is a jump--diffusion Markov process. The coupling for various values of $x$ is specified by the
two noises $W$ and $\widetilde M$ driving our SDE, which are independent of the initial condition $x$. In the case of equation \eqref{eq:CSBP}, for any sequence
$0<x_1<x_2<\cdots<x_n$, the increments $Z^{x_1},Z^{x_2}-Z^{x_1},\ldots,Z^{x_n}-Z^{x_{n-1}}$ are mutually independent. In fact this is true both concerning the increments of the processes, and the increments at some fixed value of $t$. This is the branching property. There is no reason to believe  that this independence (or equivalently, the so-called branching property) still holds  when $f$ in \eqref{eq:Zinteract} is non-linear.
However, also in this case $(Z^x)$ turns out to be a path-valued Markov process parametrized by $x$.\begin{proposition}\label{pro:Markov}
Let $\{Z^x_t,\, t,x\ge0\}$ be the solution of the collection indexed by $x$ of SDEs \eqref{eq:Zinteract}. Then
$\{Z^{x}_t,\, t\ge0\}_{x>0}$ is a $D([0,+\infty))$--valued Markov process with parameter $x$.
\end{proposition}
\bpf For $x,y>0$, let $V^{x,y}_t:=Z_t^{x+y}-Z_t^x$. It is not hard to see that $V^{x,y}$ solves the SDE
\begin{equation}\label{eq:xy}
\begin{split}
V^{x,y}_t\!&=\!y\!+\!\!\int_0^t\![f(Z^x_r+V^{x,y}_r)-f(Z^x_r)]dr\!+\!\sqrt{2\beta}\!\!\int_0^t\int_0^{V^{x,y}_r}\!\!\!\!W(dr,Z^x_r+du)\\
&\quad+\int_0^t\int_0^{V^{x,y}_{r-}}\!\!\int_0^\infty z\widetilde{M}(dr,Z^x_r+du,dz)\, ,
\end{split}
\end{equation}
where the pair of noises $({W}(dr,Z^x_r+du),\widetilde{{M}}(dr,Z^x_r+du,dz))$ is independent of $\{Z^{x'},\, 0<x'\le x\}$ and has the same law as the pair $(W,\widetilde M)$. The independence property follows from the fact that 
the restrictions of $(W,\widetilde{M})$ to disjoint sets are independent. Since the time dependent drift $v\mapsto f(Z^x_r+v)-f(Z^x_r)$
is a function of $Z^x$, and the noise terms are functions of both the solution $V^{x,y}$ and noises which are independent of $\{Z^{x'},\, 0<x'\le x\}$, we conclude that the condition law of $V^{x,y}$ given  $\{Z^{x'},\, 0<x'\le x\}$ is a function of $Z^x$. The result follows.
\epf
\section{The case without interaction}\label{sec:csbp}
Our starting point in this section will be the case $f(x)=-\alpha x$, $\alpha\ge0$ in \eqref{eq:Zinteract}, with a CSBP $Z^x$ solving \eqref{eq:CSBP}, and the corresponding L\'evy process $X$. First we recall some basic facts about the latter.
\subsection{The L\'evy process $X$}\label{sec1.2}
The branching mechanism of the CSBP $Z^x$ solving \eqref{eq:CSBP} is given as
\begin{equation}\label{psi2}
\psi(\lambda)=\alpha \lambda+\beta\lambda^2+\int_0^\infty(e^{-\lambda z}-1+\lambda z)\pi(dz).
\end{equation}

 The  Laplace transform of the associated L\'evy process $X$ is given as
\begin{equation}\label{psi1}
\E\left(\exp(-\lambda X_s)\right)=\exp(s\psi(\lambda)),\quad s,\lambda\ge0,
\end{equation}
with characteristic exponent $\psi=\psi_{\alpha,\beta,\pi}$ given by \eqref{psi2}. Our assumptions on $\beta$ and $\pi$ have been formulated in \eqref{hyp:pi}.

Let $B$ be a standard Brownian motion, 
$N$  be a Poisson random measure on $(0,+\infty)^2$ with mean measure $ds\, \pi(dz)$,  where $\pi$ satisfies \eqref{hyp:pi}, and let $\tilde{N}$ denote the compensated measure $\tilde{N}(dr,dz)=N(dr,dz)-dr\pi(dz)$. Then $X$ has the representation
\begin{equation}\label{X}
X_s=-\alpha s+\sqrt{2\beta} B_s+\int_0^s\int_0^\infty z\tilde{N}(dr,dz), \quad s \ge 0.
\end{equation}
For part of our results, we will assume that $X$ does not drift to $+\infty$, which in the presence of condition \eqref{hyp:pi}
is equivalent to
\begin{equation}\label{recur}
-\alpha  = \E(X_1)\le0.
\end{equation}
We note that our standing assumption $\beta>0$ implies that
\begin{equation}\label{conditHcont}
\int_1^\infty \frac{d\lambda}{\psi(\lambda)}<\infty.
\end{equation}
Indeed, since $e^{-\lambda z}-1+\lambda z\ge0$, we have $\psi_{\alpha,\beta,\pi}(\lambda)\ge\psi_{\alpha,\beta,0}(\lambda)=\alpha\lambda+\beta\lambda^2$.

Property  \eqref{conditHcont} implies continuity of the height process $H$ even in the case $\beta =0$, see Duquesne and Le Gall \cite {DLG}, Theorem 1.4.3. In particular, for the case $\beta >0$ considered in the present work, the height process $H$, which is then given by \eqref{1-4}, is continuous.

For the remainder of this section we assume  that \eqref{recur} holds, so that the L\'evy process $X$ hits $-x$ in finite time a.s., for any $x>0$.
We are now going to establish  properties of $X$ which will be essential for our representation of
the height process.   In the next statement, we shall write $\int_a^b$ to mean $\int_{(a,b]}$, except when $b=\infty$, in which case
$\int_a^\infty=\int_{(a,\infty)}$.
\begin{proposition}\label{ident}
For any $s>0$, $0\le a<b\le\infty$, we have
\begin{align*}
\E\int_0^s\int_a^b(z+\inf_{r\le u\le s}X_u-X_r)^+N(dr,dz)=
\E\int_0^s dr\int_a^b(z+\inf_{0\le u\le r}X_u)^+\pi(dz).
\end{align*}
\end{proposition}
\bpf
{\sc First step : $\pi(0,\infty)<\infty$.} In this case, we can write $N=\sum_{i=1}^\infty\delta_{(R_i,Z_i)}$, where $0<R_1<R_2<\cdots$ are stopping times.
Let $\F_s=\sigma\{X_r,\ 0\le r\le s\}$. Since $Z_i$ is $\F_{R_i}$--measurable, we have

\begin{align*}
\E\int_0^s\int_a^b(z+&\inf_{r\le u\le s}X_u-X_r)^+N(dr,dz)\\&=
\sum_{i=1}^\infty \E\left[{\bf1}_{\{R_i\le s, a<Z_i\le b\}}(Z_i+\inf_{R_i\le u\le s}X_u-X_{R_i})^+\right]\\
&=\sum_{i=1}^\infty \E\left[{\bf1}_{\{R_i\le s, a<Z_i\le b\}}\E\left\{(Z_i+\inf_{R_i\le u\le s}X_u-X_{R_i})^+\Big|\F_{R_i}\right\}\right]\\
&=\sum_{i=1}^\infty \E\left[{\bf1}_{\{R_i\le s, a<Z_i\le b\}}\E\left\{(z+\inf_{R_i\le u\le s}X_u-X_{R_i})^+\Big|\F_{R_i}\right\}\Big|_{z=Z_i}\right]\\
&=\sum_{i=1}^\infty \E\left[{\bf1}_{\{R_i\le s, a<Z_i\le b\}}\E\left\{(z+\inf_{0\le u\le s-r}X_u)^+\right\}\Big|_{r=R_i,z=Z_i}\right]\\
&=\E\int_0^s\int_a^b\E\left[(z+\inf_{0\le u\le s-r}X_u)^+\right]{ N(dr,dz)}\\
&=\int_0^sdr\int_a^b \E\left[(z+\inf_{0\le u\le s-r}X_u)^+\right]\pi(dz)\\
&=\int_0^sdr\int_a^b \E\left[(z+\inf_{0\le u\le r}X_u)^+\right]\pi(dz),
\end{align*}
where we have used the strong Markov property of $X_s$ for the 4th equality.

\noindent{\sc Second step : the general case.} This step is necessary only in the case $a=0$, which we now assume.
It follows from the first step that for any $k\ge1$,
\[\E\int_0^s\int_{1/k}^b(z+\inf_{r\le u\le s}X_u-X_r)^+N(dr,dz)=
\E\int_0^s dr\int_{1/k}^b(z+\inf_{0\le u\le r}X_u)^+\pi(dz).\]
We can take the limit in that identity as $k\to\infty$, thanks to the monotone convergence theorem
applied to the two expressions. \epf

\begin{lemma}\label{le:fluc}
For any $s,x>0$, we have, with $c=(1-e^{-1})^{-1}$
\[ \P\left(-\inf_{0\le r\le s}X_r\le x\right)\le \left( \frac{c}{\sqrt{\beta s}} x\right) \wedge 1. \]
\end{lemma}
\bpf
Let
\begin{align*}
T_x&=\inf\left\{s>0,\  \inf_{0\le r\le s}X_r<-x\right\}.
\end{align*}
Translating Theorem VII.1 from Bertoin \cite{B} written for spectrally negative L\'evy processes
into a statement for spectrally positive L\'evy processes, we deduce that $\{T_x,\ x\ge0\}$ is a subordinator
with the Laplace transform
\[ \E e^{-\lambda T_x}=e^{-x\Phi(\lambda)},\]
where $\Phi=\psi^{-1}$ is the inverse of the Laplace exponent $\psi$.

Combining the Markov inequality applied to the increasing function $y\to 1-e^{-y}$
and the inequality $1-e^{-y}\le y$, we get
\begin{align*}
\P\left(-\inf_{0\le r\le s}X_r\le x\right)&=\P(T_x>s)\\
&\le (1-e^{-1})^{-1}\E\left(1-e^{-T_x/s}\right)\\
&=(1-e^{-1})^{-1}\left(1-e^{-x\Phi(1/s)}\right)\\
&\le (1-e^{-1})^{-1}\Phi(1/s) x.
\end{align*}
As we have already noted, $\psi_{\alpha,\beta,\pi}(\lambda)\ge\psi_{\alpha,\beta,0}(\lambda)=\alpha\lambda+\beta\lambda^2\ge\beta\lambda^2$ since $\alpha\ge0$ (see our assumption \eqref{recur}).
 Consequently $\Phi(u)\le\sqrt{u/\beta}$ and $\Phi(1/s)\le(\beta s)^{-1/2}$. The result follows.
\epf

\begin{proposition}\label{basic:estim}
For any $s,z>0$, we have with the constant $c$ from Lemma~\ref{le:fluc}
\[ \E\left[(z+\inf_{0\le r\le s}X_r)^+\right]\le \left(\frac{c}{2\sqrt{\beta s}} z^2\right)\wedge z.\]
\end{proposition}
\bpf
It is plain that
\begin{align*}
 \E\left[(z+\inf_{0\le r\le s}X_r)^+\right]&=\int_0^z\P\left(z+\inf_{0\le r\le s}X_r\ge x\right)dx\\
 &=\int_0^z\P\left(-\inf_{0\le r\le s}X_r\le z-x\right)dx\\
 &=\int_0^z \P\left(-\inf_{0\le r\le s}X_r\le x\right)dx.
 \end{align*}
 The result now follows from Lemma \ref{le:fluc}. \epf

Next we establish the
\begin{proposition}\label{finite}
Under condition \eqref{recur}, for any $s>0$ and $0\le a<b\le\infty$,
\begin{align*} \E\int_0^sdr\int_a^b(z+\inf_{0\le u\le r}X_u)^+\pi(dz) \le C(s) \int_a^b(z\wedge z^2)\pi(dz),
\end{align*}
with $C(s)=(c\sqrt{s/\beta})\vee s$ and $c=e/(e-1)$.
\end{proposition}
\bpf
We deduce from Proposition \ref{basic:estim} and Fubini's Theorem
that
\begin{align*}
 \E\int_0^sdr\int_a^b(z+\inf_{0\le u\le r}X_u)^+\pi(dz)\le\int_a^b\pi(dz)\int_0^s\left(\frac{c}{2\sqrt{\beta r}} z^2\right)\wedge z\, dr,
  \end{align*}
from which the result follows. \epf

We now deduce readily from Propositions \ref{ident} and \ref{finite}
\begin{corollary}\label{cor:estim}
For any $s>0$ and $0\le a<b\le\infty$ we have, with $C(s)$ as in Proposition \ref{finite}, 
\begin{align*}
\E\int_0^s\int_a^b(z+\inf_{r\le u\le s}X_u-X_r)^+N(dr,dz)\le  C(s)\int_a^b(z\wedge z^2)\pi(dz).
\end{align*}
\end{corollary}


 \begin{remark}\label{theH}
By Proposition \ref{ident} and \ref{finite}, the process
\[s\mapsto U_s:=\int_0^s\int_0^\infty(z+\inf_{r\le u\le s}X_u-X_r)^+N(dr,dz)\]
is well defined. In particular, if $\pi$ is a finite measure, the 
process $U$ has only finitely many jumps on each bounded interval, so it has a
right--continuous modification. We shall from now on only consider such a
modification. In the general case \eqref{hyp:pi}, the existence of a right--continuous modification will follow from
the fact that $X$ is right--continuous and $H$ is continuous, see Proposition~\ref{convH}.

Note however that, if the measure $\pi$ obeys $\int_0^1z\, \pi(dz)=\infty$, then  the  process $U$ has infinite variation. Indeed, the contribution of the total variation  of $U$ on the interval 
$[r,s]$ induced by a jump of size $z$ of $X$ at some time $r'\in(r,s)$ is bounded from below by $z$ and from above by $2z$. Consequently the total variation $TV_U([r,s])$ of $U$ on the interval $[r,s]$ satisfies
\[ \int_r^s\int_0^\infty z N(dr,dz)\le TV_U([r,s])\le 2   \int_r^s\int_0^\infty z N(dr,dz).\]
It follows from well--known properties of Poisson random measures that $\int_r^s\int_0^\infty z N(dr,dz)=+\infty$
a.s., unless $\int_0^\infty (z\wedge 1) \pi(dz)<\infty$. 
\end{remark}

\subsection{The  case $\pi=0$}
In this subsection we assume that the L\'evy process $X$ is continuous, i.e.
\[ X_s=-\alpha s+\sqrt{2\beta} B_s, \quad  s\ge 0. \]
\begin{proposition}\label{pro:Hcont}
In the case $\pi=0$, we have
\[H_s=\frac{1}{\beta}\left(X_s-\inf_{0\le r\le s}X_r\right), \quad  s\ge 0.\]
\end{proposition}
\bpf This result follows readily from \eqref{1-4}, since 
$r\to \overline{X}^s_r$ is continuous and increases from
$\inf_{0\le r\le s}X_r$ to $X_s$. \epf

In this case, $H$ is a drifted Brownian motion reflected above $0$, and thus a fortiori a continuous semi-martingale. The next proposition states the second Ray--Knight theorem for this particular case. Let us define $L^t(s)=L^t(s,H)$ and
\begin{equation}\label{defSx}
S_x=\inf\{s>0,\ L^0(s,H)>x\}.
\end{equation}
\begin{proposition}\label{RN-cont}
The process $\{L^t(S_x),\ t\ge0\}$ is a CSBP with branching mechanism $\psi_{\alpha,\beta,0}$, starting from $x$ at time $t=0$.
\end{proposition}
\bpf
This  is classical, see e.g. Revuz and Yor \cite{RY} Chapter XI \S 2, and   
Theorem 5.1 in Ba, Pardoux, Sow \cite{BPS} for an identification of the constants in our case.
\epf

\begin{remark} We note that the scaling of the local time of $H$ is such that 
$L^t(s,H)=\frac{\beta}{2}\mathcal{L}^t(s,H)$, where $\mathcal{L}^t(s,H)$ is the semi--martingale local time as defined in Revuz and Yor \cite{RY} (see Corollary VI.1.9, page 227). 
Then, from the Tanaka formula, see Theorem VI.1.2 page 222 in \cite{RY}, and Proposition~\ref{pro:Hcont}, we have
\begin{align*}
 H_s&=\int_0^s{\bf1}_{H_r>0}dH_r+\frac{1}{\beta}L^0(s,H)\\
 &=\frac{1}{\beta}X_s+\frac{1}{\beta}L^0(s,H).
 \end{align*}
  The second equality can be justified as follows.  Proposition \ref{pro:Hcont} tells that
 \begin{align*}
 \beta \, dH_s=dX_s+d(-\inf_{r\le s}X_r).
 \end{align*}
 However it is plain that 
 \[ {\bf1}_{H_s>0}\, d(-\inf_{r\le s}X_r)=0,\]
 since $\inf_{r\le s}X_r$ decreases only when $H_s=0$. Consequently
 \begin{align*}
 \beta \int_0^s{\bf1}_{H_r>0}dH_r =\int_0^s{\bf1}_{H_r>0}dX_r =X_s,
 \end{align*}
since ${\bf1}_{H_r>0}=1$ for Lebesgue-a.a. $r$ and $X$ is a drifted Brownian motion. We note  in particular that $L^0(s,H)=-\inf_{0\le r\le s}X_r$, which is L\'evy's correspondence between the local time of a reflected BM at the origin and the current minimum of  a~BM.
 \end{remark}

\subsection{The case of finite $\pi$}
We now suppose that $\pi$ is a finite measure. In that case, in view of condition~\eqref{hyp:pi}, $z\pi(dz)$ is also a finite measure,
and if we let
\[  \alpha'=\alpha+\int_0^\infty z\pi(dz),\]
we have that
\begin{equation}\label{eq:x}
 X_s=\sqrt{2\beta} B_s+P_s-\alpha's,
 \end{equation}
 where
\begin{equation}\label{eq:cpp}
 P_s=\int_0^s\int_0^\infty zN(dr,dz), \quad s\ge 0,
 \end{equation}
is a compound Poisson process.

Recall the notation introduced in \eqref{1-4}.
We note that $[0,s]\ni r\mapsto  \overline{X}^s_r$ is increasing. Denote by $\Delta \overline{X}^s_r$ its possible jump at time $r$. It follows readily  
from \eqref{1-4} that
\begin{equation}\label{1.4}
\beta H_s=X_s- \overline{X}^s_0-\sum_{0\le r\le s}\Delta \overline{X}^s_r\, ,
\end{equation}
which we rewrite as
\begin{equation}\label{betaH}
\beta H_s=X_s-\inf_{0\le r\le s}X_r-\int_0^s\int_0^\infty\left(z+\inf_{r\le u\le s}(X_u-X_r)\right)^+N(dr,dz),
\end{equation}
hereby using the equality
\begin{equation}\label{identif}
 \int_0^s\int_0^\infty\left(z+\inf_{r\le u\le s}(X_u-X_r)\right)^+N(dr,dz)=\sum_{0\le r\le s}\Delta \overline{X}^s_r.
 \end{equation}
As observed in Remark \ref{theH}, the third term on the right-hand
side of \eqref{betaH} has only finitely many jumps on each bounded
interval and its jumps compensate those of the process $X$.
This shows that $H$ is continuous. The fact that the last term in \eqref{betaH} has bounded variation shows that $H$ is a semi--martingale.
We have thus proved
\begin{proposition}\label{heightcpp}
If the L\'evy process $X$ is given by  \eqref{eq:x} and \eqref{eq:cpp} with $\pi$ finite,
then the associated height process $H$ is given by \eqref{betaH}, and it
is a continuous semi--martingale.
\end{proposition}
Recall that the
second term on the right of \eqref{betaH} reflects the process above
0. We will explain in words what the third term in \eqref{betaH} does. For
that purpose, we need to define for $0\le r\le s$
\[ \tilde{X}^r_s=X_s-\inf_{0\le r\le s}X_r-\int_0^r\int_0^\infty\left(z+\inf_{u\le v\le r}(X_v-X_u)\right)^+N(du,dy), \]
which is the same as $\beta H_s$, except that we have stopped the third term at time $r$.
At each jump time $r$ of $X$, draw a piece of horizontal
line which starts from $(r,\beta H_r)$, and extends to $r^+:=\inf\{s>r, \tilde{X}^r_s\le\beta H_r\}$. Then pull down
$\tilde{X}^{r}_r$ to $\tilde{X}^{r}_{r-}$, and reflect the piece of trajectory of $\{\tilde{X}^r_u, r\le u\le r^+\}$
above the level of $\tilde{X}^{r}_{r-}$, that is above the ``horizontal stick'' which extends from 
$(r,\tilde{X}^{r}_{r-})$
to $(r+,\tilde{X}^{r}_{r+})$.

We are now going to give a new derivation of the Ray--Knight theorem in this case, since our proof
of the corresponding result in the case with interaction will be based upon the same argument.

Let $X$ be the L\'evy process given by \eqref{psi1}, started in $0$ and stopped at the time $S_x$ when first hitting $-x$. Let the height process $H$ of $X$ be given by \eqref{1-4} (or equivalently by \eqref{1.4}), and $L^t(s)$ be the local time accumulated by $H$ at height $t$ between times $0$ and $s$.

\begin{proposition}\label{RNcpp}
The process $\{L^t(S_x),\ t\ge0\}$ is a CSBP with branching mechanism $\psi:= \psi_{\alpha,\beta,\pi}$, starting from $x$ at time $t=0$.
\end{proposition}
\bpf
Let us first recall how we can construct $X$ and $S_x$ iteratively from (pieces of) drifted Brownian motions together with (atoms of) the Poisson process $\Pi$ with intensity $ds \,  \pi(dz)$. For $y >0$ let $B^y$ be a  BM with drift $\alpha -\int_0^\infty z\pi(dz)$ started in $y$ and stopped when first hitting $0$; let us denote this hitting time by $\mathcal S(B^y)$.

 The first 
step in the iteration is $X^{(0)}:= -x+B^x$. Let $S_x^{(0)}$ be the time at which $X^{(0)}$ first hits $-x$; note that $S_x^{(0)}=\mathcal S(B^x)$. Consider a Poisson process $\Pi_0$ on $[0,S_x^{(0)}]\times(0,+\infty)$ with intensity $ds \,  \pi(dz)$. Denote the points of $\Pi_0$ by $(s_i, z_i)_{1\le i\le J}$. If $J=0$, the iteration stops at step $0$. Otherwise each atom $(s_i, z_i)$ gives rise to the injection of a path $B^{z_i}$ (defined on an interval of length $\mathcal S(B^{z_i})$) into $X^{(0)}$ as described below for $ B^{z_m}$. Each piece $B^y$ that is injected in the $k$-th iteration is defined on some interval $I$, and gives rise to a Poisson process on $I\times \mathbb R_+$, whose points in turn give rise to new injected pieces. This procedure terminates after finitely many steps, ending in $X$.

Let $m \in \{1,\ldots, J\}$ be such that $H^{(0)}_{s_m} = \min\{H^{(0)}_{s_i}: 1\le i \le J\}$. Given that $\mathcal S(B^{z_m}) =\tilde s$, take a Poisson process  $\tilde \Pi_1$ with intensity measure $ds\, \pi(dz)$ on $[s_m, s_m+\tilde s]\times\mathbb R_+$. Then transport $\Pi_0$ into  $\tilde \Pi_0$ by keeping each point $(s_i,z_i)$ with $s_i < s_m$ as it is, and  shifting each point $(s_i,z_i)$ with $s_i > s_m$ into $(s_i+\tilde s, z_i)$. Put $\Pi_1:= \tilde \Pi_0 + \tilde \Pi_1$, and keep iterating. 

Let  $H^{(0)}$ be the height process of $X^{(0)}$,  given by \eqref{1.4} with $X^{(0)}$ instead of $X$. 
 Inject  $B^{z_m}$ into $X^{(0)}$, by defining 
\[ X^{(1)}_s=\begin{cases} X^{(0)}_s&, \text{for } 0\le s\le s_m,\\
                 X_{s_m}^{(0)}+ B^{z_m}_{s-s_m}&, \text{for } s_m\le s\le s_m+\mathcal  S(B^{z_m}),\\
                 X_{s-\mathcal  S(B^{z_m})}^{(0)}&, \text{for } s_m+\mathcal  S(B^{z_m}) < s\le S^{(0)}_x+\mathcal  S(B^{z_m}) =:S^{(1)}_x .
                 \end{cases}
                 \]
Let  $H^{(1)}$ be the height process of $X^{(1)}$,  given by \eqref{1.4} with $X^{(1)}$ instead of~$X$. We note that $H^{(1)}_{s_m}=H^{(1)}_{s_m +\mathcal  S(B^{z_m})} = H^{(0)}_{s_m}$, and that (with $T_1:=H^{(0)}_{s_m}$) we have $L^{T_1}(H^{(1)}, S^{(1)}_x) = L^{T_1}(H^{(0)}, S^{(0)}_x)+z_m$.

A key observation is that the reflection of $H$ below $T_1$ equals the reflection of $H^{(0)}$ below $T_1$, and that $L^t(H,S_x)=L^t(H^{(0)},S^{(0)}_x)$ for $0\le t<T_1$, $L^{T_1}(H,S_x)= L^{T_1}(H^{(1)}, S^{(1)}_x)$. Consequently,
on $\{t < T_1\}$ we have $L^{t}(H,S_x)= L^{t}(H^{(0)}, S^{(0)}_x) $. The height of the lowest jump of the local time of $H$ is $T_1$, which is measurable with respect to $(X^{(0)},  \Pi_0)$. By the classical Ray-Knight theorem (Proposition \ref{RN-cont}), $L^t(H,S_x)$ follows before its first jump the dynamics of a subcritical Feller branching diffusion. Moreover,
\begin{align*}
\P(T_1 > t | X^{(0)}) &= 
\exp\left(- \pi(\mathbb R_+)\int_0^{S^{(0)}_x} 1_{\{H_s^{(0)} \le t \}}ds\right)\\
&= \exp\left(- \pi(\mathbb R_+)\int_0^t L^v(H^{(0)},S^{(0)}_x)dv\right),
\end{align*}
which shows that the first jump $T_1$ of $\{L^t(S_x),\ t\ge0\}$ comes at  rate $\pi(\mathbb R_+) L^t(S_x) dt$, since
$\int_0^t L^v(H^{(0)},S^{(0)}_x)dv=\int_0^t L^v(H,S_x)dv$ when $t<T_1$. Also, its size has distribution $\pi /\pi(\mathbb R_+)$. Thus, up to and including $T_1$, $\{L^t(S_x),\ t\ge0\}$ is a CSBP with branching mechanism $\psi$.  Proceeding in the same manner from $T_1 = H^{(1)}_{s_m}$ upwards, we arrive at our assertion. 
 \epf

\subsection{The  general case}\label{sec:gen_case}
\subsubsection{The height process}\label{Sec241}
We now consider the general case, that is $\pi$ satisfies \eqref{hyp:pi}. Consequently, for any $\eps>0$, $\pi(\eps,\infty)<\infty$.
We  define $\pi_k(dz)={\bf1}_{(\eps_k,\infty)}(z)\pi(dz)$, where $\eps_k$ is a sequence of positive reals which decreases to $0$, and
\[  \psi_k=\psi_{\alpha,\beta,\pi_k}.   \]
The corresponding L\'evy process $X^k$ admits the L\'evy--It\^o decomposition
\begin{align*}
X^k_s&=-\alpha s+\sqrt{2\beta}B_s+\int_0^s\int_{\eps_k}^\infty z\tilde{N}(dr,dz)\\
&=-(\alpha+\int_{\eps_k}^\infty z\pi(dz))s+\sqrt{2\beta}B_s+\int_0^s\int_{\eps_k}^\infty z N(dr,dz).
\end{align*}
The last term in the right--hand side is  a compound Poisson process.
We have
\begin{lemma}\label{convX}
As $k\to\infty$, $X^k_s\to X_s$ in $L^1(\Omega)$, locally uniformly with respect to $s$.
\end{lemma}
\bpf
It is plain that
\begin{align*}
\E\left[\sup_{0\le r\le s}|X_r-X^k_r|\right]&\le\left(\E\left[\sup_{0\le r\le s} \left|\int_0^r\int_0^{\eps_k}z\tilde{N}(du,dz)\right|^2 \right]\right)^{1/2}\\
&\le 2\, \sqrt{s \int_0^{\eps_k}z^2\pi(dz)}\\
&\to0,
\end{align*}
as $k\to\infty$, where we have used Doob's inequality. The result follows.  \epf

Thanks to Proposition \ref{heightcpp}, the height process $H^k$ associated to the L\'evy process $X^k$ is given by
\begin{equation}\label{eq:Hk}
\beta H^k_s=X^k_s-\inf_{0\le r\le s}X^k_r-\int_0^s\int_{\eps_k}^\infty(z+\inf_{r\le u\le s}X^k_u-X^k_r)^+N(dr,dz).
\end{equation}
Under our standing assumption \eqref{hyp:pi} we have
\begin{proposition}\label{convH}
For any $s\ge0$, $H^k_s\to H_s$ in
probability, where  
$H$ is given by either of the   formulas \eqref{betaH}, \eqref{1.4} or \eqref{1-4}, and is continuous.
\end{proposition}
\bpf
Lemma \ref{convX} implies that $X^k_s-\inf_{0\le r\le s}X^k_r\to X_s-\inf_{0\le r\le s}X_r$ in probability,
locally uniformly in $s$. We now consider the last term in \eqref{eq:Hk} and prove pointwise convergence.
It follows from Corollary \ref{cor:estim} with $a=0$ and $b=\eps_k$ that
\begin{align*}
\E\int_0^s\int_0^{\eps_k}(z+\inf_{r\le u\le s}X_u-X_r)^+N(dr,dz)\le  C(s)
\int_0^{\eps_k}(z\wedge z^2)\pi(dz),
\end{align*}
which clearly tends to 0, as $k\to\infty$. From an adaptation of the argument of Proposition \ref{ident}, we deduce that
\begin{align*}
\E&\int_0^s\int_{\eps_k}^\infty\left|(z+\inf_{r\le u\le s}(X_u-X_r))^+ -(z+\inf_{r\le u\le s}(X^k_u-X^k_r))^+\right|N(dr,dz)\\
&=\E\int_0^sdr\int_{\eps_k}^\infty\left|(z+\inf_{0\le u\le r}X_u)^+ -(z+\inf_{0\le u\le r}X^k_u)^+\right|\pi(dz)\\
&\le\E\int_0^sdr\int_{\eps_k}^\infty\Bigg\{\left(z+\inf_{0\le u\le r}X_u\vee\inf_{0\le u\le r}X^k_u\right)^+
\\ &\phantom{AAAAAAAAAAAAAAAAAAAAAAA} \wedge \left|\inf_{0\le u\le r}X_u-\inf_{0\le u\le r}X^k_u\right|\Bigg\}
\pi(dz),
\end{align*}
hence
\begin{equation}\label{estim1}
\begin{split}
\E&\int_0^s\int_{\eps_k}^\infty\left|(z+\inf_{r\le u\le s}(X_u-X_r))^+ -(z+\inf_{r\le u\le s}(X^k_u-X^k_r))^+\right|N(dr,dz)\\
&\le\int_0^sdr\int_{\eps_k}^\infty\E\left[\left(z+\inf_{0\le u\le r}X_u\vee\inf_{0\le u\le r}X^k_u\right)^+\right]
\\& \phantom{AAAAAAAAAAAAAAAAAAAAAAA} \wedge \E\left|\inf_{0\le u\le r}X_u-\inf_{0\le u\le r}X^k_u\right|
\pi(dz).
\end{split}
\end{equation}
We deduce from Lemma \ref{convX} that
\begin{equation}\label{eq:convX}
\E\, \left|\inf_{0\le u\le r}X_u-\inf_{0\le u\le r}X^k_u\right|\to0,
\end{equation}
as $k\to\infty$.
Arguing as in the proof of Proposition \ref{basic:estim}, we obtain
\begin{align*}
\E&\left[\left(z+\inf_{0\le u\le r}X_u\vee\inf_{0\le u\le r}X^k_u\right)^+\right]\\&\qquad=\int_0^z\P\left(-\left[\inf_{0\le u\le r}X_u\vee\inf_{0\le u\le r}X^k_u\right]\le x\right)dx\\
&\qquad\le\int_0^z\left[\P\left(-\inf_{0\le u\le r}X_u\le x\right)+\P\left(-\inf_{0\le u\le r}X^k_u\le x\right)\right]dx\\
&\qquad\le c\left[\Phi\left(\frac{1}{r}\right)+\Phi_k\left(\frac{1}{r}\right)\right]z^2  \, \le \, \frac{c}{\sqrt{\beta s}}  z^2.
\end{align*}
It is plain that the left hand side in  the previous chain of inequalities  is dominated by $z$, hence we have proved that
\begin{equation}\label{ineq:inf}
\E\left[\left(z+\inf_{0\le u\le r}X_u\vee\inf_{0\le u\le r}X^k_u\right)^+\right]\le\left( \frac{c}{\sqrt{\beta r}}  z^2\right)\wedge z.
\end{equation}
The right--hand side of \eqref{ineq:inf} is $dr\times\pi(dz)$\,-\,integrable
over $[0,s]\times(0,\infty)$ for any $s>0$. It then follows from
\eqref{eq:convX} and the dominated convergence theorem that the left--hand
side of \eqref{estim1} tends to 0 as $k\to\infty$. We can now take the
limit in \eqref{eq:Hk}, yielding the convergence. 

It is clear that \eqref{identif} still holds in the general situation, which re--establishes the  formulas \eqref{1.4} and \eqref{1-4}. From \eqref{1.4} the continuity of $H$ is essentially clear, as claimed in \cite{DLG}.
Let us give a quick explanation. The right continuity follows from the right continuity of the three terms on the right of \eqref{1.4}. The left continuity follows from that of the second term, while the eventual jumps of the first and the third term compensate.
\epf

Note that, under condition \eqref{conditHcont} which is weaker than $\beta>0$, Duquesne and Le Gall \cite{DLG}, Sec.1.4.3, prove that $H$ is H\"older continuous. We shall not need that property.

We first prove
\begin{lemma}\label{le:incH}
For any $\bar{s}>0$,
there exists a random increasing function $\Psi:[0,1]\mapsto \R_+$ such that $\Psi(h)\downarrow0$ a.s. as $h\downarrow0$, and for any $0\le s\le \bar{s}$, any $0< h\le1$,
\[ (H^k_{s+h}-H^k_s)_-\le\Psi(h),\,  \forall k\ge1.\]
\end{lemma}
\bpf
Since $X$ is a L\'evy process with only positive jumps, it is not hard to check by contradiction that  
\[ \Phi_X(h):=\sup_{0\le r\le \bar{s}, 0\le s-r\le h}(X_s-X_r)_-\]
is a.s. a continuous function of $h$ on $[0,1]$ such that $\Phi_X(0)=0$.
Since $X^k\to X$ uniformly in probability on $[0,\bar{s}]$, one obtains that $\Phi_{X^k}(h)\to\Phi_X(h)$ in probability as $k\to\infty$, for any $h>0$. Since each $\Phi_{X^k}$ is increasing and the limit is continuous, it follows from the second Dini theorem, that the convergence  in probability is uniform  w.r.t. $h\in[0,1]$, see
 the statement 127 on page 81, and the proof on page 270 in P\'olya and Szeg\"o \cite{PS}. This implies readily that
\[ \Psi(h):=\beta^{-1}\left(\sup_{k\ge1}\Phi_{X^k}(h)\vee\Phi_X(h)\right), \quad h \ge 0\]
is a.s. continuous in $h$, and $\Psi(0)=0$.

Now, using \eqref{1.4} and abbreviating $Y:= X^k$
\begin{align*}
\beta(H^k_{s+h}-H^k_s)&=Y_{s+h}-Y_s-\overline{Y}^{s+h}_0+\overline{Y}^s_0
-\sum_{0\le r\le s}(\Delta \overline{Y}^{s+h}_r-\Delta\overline{Y}^s_r)\\
&\quad-\sum_{s<r\le s+h}\Delta\overline{Y}^{s+h}_r\\
&\ge Y_{s+h}-Y_s-\sum_{s<r\le s+h}\Delta\overline{Y}^{s+h}_r.
\end{align*}
Since $Y_{s+h}-\sum_{s<r\le s+h}\Delta\overline{Y}^{s+h}_r\ge \inf_{s\le r\le s+h}Y_r$, we conclude that
\[ \beta(H^k_{s+h}-H^k_s)\ge\inf_{s\le r\le s+h}Y_r-Y_s,\] 
and consequently
\[(H^k_{s+h}-H^k_s)_-\le \Psi(h),\]
which proves the result.
\epf

We now deduce from the two previous statements
\begin{corollary}\label{cor:convH}
Under the above assumptions, $H^k_s\to H_s$ in probability, locally uniformly w.r.t. $s$.
\end{corollary}
Corollary \ref{cor:convH} is an immediate consequence of Proposition \ref{convH}, Lemma~\ref{le:incH}, the following extension of the second Dini theorem, and the equivalence of convergence in probability and the fact that from any subsequence, one can extract a further subsequence which converges a.s..
\begin{lemma}
Consider a sequence $\{g_k,\ k\ge1\}$ of functions from $\R_+$ into~$\R$ and $T>0$, which are such that for any $0\le t\le T$,
$g_k(t)\to g(t)$, where $g:[0,T]\mapsto\R$ is continuous, and
$\sup_{k\ge1}(g_k(t+h)-g_k(t))_- \to0$, as $h\to0$. Then $g_k(t)\to g(t)$ uniformly w.r.t. $t\in[0,T]$.
\end{lemma}
\bpf  Let 
 $\eps>0$ be arbitrary. 
 Since $t\mapsto g(t)$ is uniformly continuous on the compact interval
$[0,T]$, there exists 
$\eta>0$ small enough such that whenever $s,t\in[0,T],\ 0<t<s<t+\eta$, 
\begin{align}
|g(s)-g(t)|&\le \frac{\eps}{3}\label{contH},\\
g_k(s)-g_k(t)&\ge-\frac{\eps}{3}\ \text{ for all }k\ge1,\label{incH}
\end{align}
where the second inequality follows from our assumption on the sequence $g_k$.

We next choose an integer $N>T/\eta$ and $0=t_0<t_1<\cdots<t_N=T$ such that $t_{j+1}-t_j<\eta$, for all $0\le j<N-1$.
We now choose $k_\eps$ large enough such that, for any $k\ge k_\eps$, $1\le j\le N$,
\begin{equation}\label{convHkH}
|g_k(t_j)-g(t_j)|\le\frac{\eps}{3}.
\end{equation}
Now for any $t\in[0,T]$ either $t=t_j$ for some $1\le j\le N$ (and then \eqref{convHkH} ensures that $|g_k(t)-g(t)|<\varepsilon$), or else there exists $0\le j<N$ such that
$t_j<t<t_{j+1}$. In that case we obtain, using successively  \eqref{incH}, \eqref{convHkH} and \eqref{contH},
the two following inequalities:
\begin{align*}
(i)\qquad g_k(t)&\le g_k(t_{j+1})+\frac{\eps}{3}\\
&\le g(t_{j+1})+\frac{2\eps}{3} \, \le\,  g(t)+\eps,\\
(ii) \qquad g_k(t)&\ge g_k(t_{j})-\frac{\eps}{3}\\
&\ge\,  g(t_{j})-\frac{2\eps}{3} \quad \ge \,   g(t)-\eps.
\end{align*}
The result clearly follows from those inequalities and the fact that $\eps>0$ is abitrary.
\epf

\subsubsection{The local time of the height process}
Let $L^t(s)$ denote the local time accumulated by the process $H$, defined by \eqref{1.4} or \eqref{betaH}, at level $t$ up to time $s$. The existence of $L^t(s)$ was established already in Duquesne and Le Gall \cite{DLG}.  We shall give an independent definition of  $L^t(s)$, via an It\^o--Tanaka formula for $(H-t)^+$, and prove some regularity.
%
\begin{proposition}\label{tanaka}
We have
\begin{align*}
\beta(H_s-t)^+\!=\! \int_0^s\!\!{\bf1}_{H_r>t}dX_r \!
- \! \int_0^s\!\int_0^\infty \!\!{\bf1}_{H_r>t}(z+\inf_{r\le u\le s}X_u-X_r)^+N(dr,dz)\!
+\! L^t(s),
\end{align*}
where $L^t(s)$ is for any $s>0$, $t\ge0$ the local time accumulated by $H$ at level $t$ up to time $s$, in the sense that it satisfies the occupation times formula.
\end{proposition}
The formula in the Proposition can be rewritten as
\begin{equation}\label{defLT}
L^t(s)
=\beta(H_s-t)^+\!\!-\!\!\int_0^s{\bf1}_{H_r>t}dX_r
+\int_0^s\int_0^\infty\!\!\!\!\!{\bf1}_{H_r>t}(z\!+\!\!\inf_{r\le u\le s}X_u-X_r)^+N(dr,dz).
\end{equation}
The proof of Proposition\ref{tanaka} will be based on a limiting procedure along the sequence $X^k$ of L\'evy processes  associated to 
$\pi_k(dz)={\bf1}_{z>\eps_k}\pi(dz)$. This gives us a construction of  the local time that is different from the construction in  Duquesne and Le Gall \cite{DLG}, but leads to  the same result,
as a consequence of the occupation time formula. 

Note that the corresponding height process $H^k$ is a continuous semi--martingale, whose local time 
is well-defined using the classical theory, see e.g. Chapter~VI in 
Revuz and Yor \cite{RY}. 
We have the formula, analogous to \eqref{defLT}
\begin{equation}\label{LTk}
\begin{split}
L^{t}_k(s)&=\beta(H^k_s-t)^+-\int_0^s{\bf1}_{H^k_r>t}dX^k_r\\
&\quad+\int_0^s\int_{\eps_k}^\infty\!\!\!\!{\bf1}_{H^k_r>t}(z+\inf_{r\le u\le s}X^k_u-X^k_r)^+N(dr,dz).
 \end{split}
\end{equation}
Note that the formula would be different if $L^t_k(s)$ were the ``semi--martingale local time'', as defined in
\cite{RY}. In that case there would be a factor $\frac{\beta}{2}$ in front of the local time. Indeed, after the division 
of the whole formula by $\beta$, we should find a factor $\frac{1}{2}$ in front of the local time, see 
the second formula in Theorem VI.1.2 in \cite{RY}.

Before proving the above Proposition, let us establish a technical Lemma. 
\begin{lemma}\label{le:estimLT}
For any $s>0$,
\[ \sup_{t>0,\ k\ge1}\E L^t_k(s)<\infty.\]
\end{lemma}
\bpf We need to show successively
\begin{align*}
&\sup_{t>0,\ k\ge1}\E (H^k_s-t)^+<\infty,\\
&\sup_{t>0,\ k\ge1}\E\left|\int_0^s{\bf1}_{H^k_r>t}dX^k_r\right|<\infty,\\
&\sup_{t>0,\ k\ge1}\E\int_0^s\int_{\eps_k}^\infty\!\!\!\!{\bf1}_{H^k_r>t}(z+\inf_{r\le u\le s}X^k_u-X^k_r)^+N(dr,dz)<\infty.
\end{align*}
The first estimate is an easy exercise which we leave to the reader. The third one follows readily from
\[\int_0^s\!\!\int_{\eps_k}^\infty\!\!\!\!{\bf1}_{H^k_r>t}(z+\inf_{r\le u\le s}X^k_u-X^k_r)^+N(dr,dz)
\le\int_0^s\!\!\int_{\eps_k}^\infty\!\!(z+\inf_{r\le u\le s}X^k_u-X^k_r)^+N(dr,dz)\]
and Proposition \ref{finite}. It remains to consider
\[\int_0^s{\bf1}_{H^k_r>t}dX^k_r=-\alpha\int_0^s{\bf1}_{H^k_r>t}dr+\sqrt{2\beta}\int_0^s{\bf1}_{H^k_r>t}dB_r
+\int_0^s{\bf1}_{H^k_r>t}\int_{\eps_k}^\infty z\tilde{N}(dr,dz).\]
The first term on the right is bounded in absolute value by $|\alpha|s$. We estimate the second term using Cauchy--Schwartz
\[ \E\left|\int_0^s{\bf1}_{H^k_r>t}dB_r\right|\le\sqrt{s}.\]
Finally
\begin{align*}
&\E\left|\int_0^s{\bf1}_{H^k_r>t}\int_{\eps_k}^\infty z\tilde{N}(dr,dz)\right|\\ &\le \E\left|\int_0^s{\bf1}_{H^k_r>t}\int_{\eps_k}^1 z\tilde{N}(dr,dz)\right|
+\E\left|\int_0^s{\bf1}_{H^k_r>t}\int_{1}^\infty z\tilde{N}(dr,dz)\right|\\
&\le \sqrt{s\int_0^1z^2\pi(dz)}+2s\int_1^\infty z\pi(dz).
\end{align*}
The result follows.\epf

We now turn to the

\noindent{\sc Proof of Proposition \ref{tanaka}} We first consider the case  $\int_0^\infty z\pi(dz)<\infty$ (which certainly applies to $\pi(dz):= \pi_k(dz)= {\bf1}_{z>\eps_k}\pi(dz) $).
Then $H_s$ is a continuous semi--martingale, and the formula of our Proposition follows from 
It\^o--Tanaka's formula (see e.g. the second identity in Theorem VI.1.2 in \cite{RY}), but with a 
different constant in front of the local time, due to our definition \eqref{loctime}. It is then crucial to note that
whenever we have a point $(r,z)$ of the Point Process $N$ such that $H_r\le t$, then until the first time $s$ for which
  $z+\inf_{r\le u\le s}X_u-X_r=0$, the 
process $u \mapsto \inf_{r\le v\le u}X_v$ decreases only when $H_u=H_r\le t$, hence the term ${\bf1}_{H_r>t}$
factorizes in the last integral.

We now take the limit along a sequence $X^k$  associated to 
$\pi_k$, thus establishing the It\^o--Tanaka formula in the general case.

From the occupation time formula, for any $g\in C([0,\infty))$ with compact support
\[ \int_0^\infty g(t)L^t_k(s)dt=\int_0^s g(H^k_r)dr. \]
Clearly $\int_0^s g(H^k_r)dr\to \int_0^s g(H_r)dr$ as $k\to\infty$. Denote by $R^t(s)$ (resp. $R^t_k(s)$) the right--hand side of \eqref{defLT} (resp. of \eqref{LTk}). The Proposition will clearly follow from
\begin{equation*}\label{eq:limit}
\sup_{t>0}\E\left|R^t(s)-R^t_k(s)\right|\to0,\  \text{ as } k\to\infty.
\end{equation*}
In other words, all we have to show is that, as $k\to\infty$,
\begin{align}
&\sup_{t>0}\E\left|(H_s-t)^+-(H^k_s-t)^+\right|\to0,\label{conv1}\\
&\sup_{t>0}\E\left|\int_0^s{\bf1}_{H_r>t}dX_r-\int_0^s{\bf1}_{H^k_r>t}dX^k_r\right|\to0,\label{conv2}\\\begin{split}
&\sup_{t>0}\E\Big|\int_0^s\int_0^\infty\!\!\!\!\!{\bf1}_{H_r>t}(z\!+\!\!\inf_{r\le u\le s}X_u-X_r)^+N(dr,dz)\quad
\\&\phantom{AAAAAAAA}-
\int_0^s\int_{\eps_k}^\infty\!\!\!\!\!{\bf1}_{H^k_r>t}(z\!+\!\!\inf_{r\le u\le s}X^k_u-X^k_r)^+N(dr,dz)\Big|\to0.
\label{conv3}
\end{split}
\end{align}

Since $x\to(x-t)^+$ is continuous, \eqref{conv1} follows from Corollary \ref{cor:convH} and the uniform integrability of the sequence $H^k_s$. We next establish \eqref{conv3}.
We argue similarly as in the proof of Proposition \ref{convH}.
\begin{align*}
\E&\Bigg|\!\int_0^s\!\!\int_0^\infty\!\!\!\!\!\!{\bf1}_{H_r>t}(z+\inf_{r\le u\le s}X_u-X_r)^+N(dr,dz)\\  &\phantom{AAAAA}-
\int_0^s\!\!\int_{\eps_k}^\infty\!\!\!\!\!\!{\bf1}_{H^k_r>t}(z+\inf_{r\le u\le s}X^k_u-X^k_r)^+N(dr,dz)\Bigg| \\
&\le\E\!\int_0^s\int_0^{\eps_k}\!\!\!\!\!\!\!{\bf1}_{H_r>t}(z+\inf_{r\le u\le s}X_u-X_r)^+N(dr,dz)\\
&\phantom{AAAAA}+\E\int_0^s\int_{\eps_k}^\infty\Big|{\bf1}_{H_r>t}(z+\inf_{r\le u\le s}(X_u-X_r))^+ \\ &\phantom{AAAAAAAAAAAA}-{\bf1}_{H^k_r>t}(z+\inf_{r\le u\le s}(X^k_u-X^k_r))^+\Big|N(dr,dz).
\end{align*}
The first term on the right is bounded from above by the same term without the factor ${\bf1}_{H_r>t}$, which tends to $0$ as $k\to\infty$ thanks to Corollary \ref{cor:estim}. We now estimate the second term.
\begin{align*}
\E&\int_0^s\int_{\eps_k}^\infty\left|{\bf1}_{H_r>t}(z+\inf_{r\le u\le s}(X_u-X_r))^+ -{\bf1}_{H^k_r>t}(z+\inf_{r\le u\le s}(X^k_u-X^k_r))^+\right|N(dr,dz)\\
&\le \E\int_0^s\int_{\eps_k}^\infty\left|{\bf1}_{H_r>t}-{\bf1}_{H^k_r>t}\right|(z+\inf_{r\le u\le s}(X_u-X_r))^+ N(dr,dz)\\
&\quad+\E\int_0^s\int_{\eps_k}^\infty\left|(z+\inf_{r\le u\le s}(X_u-X_r))^+ -(z+\inf_{r\le u\le s}(X^k_u-X^k_r))^+\right|N(dr,dz)
\end{align*}
The first term on the right hand side of the last inequality tends to zero by dominated convergence, while the convergence to zero of the second term was proved in Proposition \ref{convH}. In order to finally establish \eqref{conv2}, we first note that
\begin{align*}
\int_0^s{\bf1}_{H_r>t}dX_r\!-\!\int_0^s{\bf1}_{H^k_r>t}dX^k_r\!=\!
\int_0^s\!\left[{\bf1}_{H_r>t} - {\bf1}_{H^k_r>t}\right]\!dX_r
\!+\!\int_0^s{\bf1}_{H^k_r>t}\int_0^{\eps_k}\!\!\!z\tilde{N}(ds,dz).
\end{align*}
The sup in $t$ of the expectation of the second term on the right tends to zero since
\begin{align*}
\E\left[\left(\int_0^s{\bf1}_{H^k_r>t}\int_0^{\eps_k}z\tilde{N}(ds,dz)\right)^2\right]
&\le s\int_0^{\eps_k}z^2\pi(dz)\\
&\to 0,\ \text{ as } k\to\infty.
\end{align*}
Concerning the first term, all we need to do is to use the same decomposition and the same kind of estimates as used in the proof of Lemma \ref{le:estimLT}, combined with the following 
\begin{equation}\label{eq:conv:0}
\sup_{t>0}\E\int_0^s\left|{\bf1}_{H_r>t} - {\bf1}_{H^k_r>t}\right|dr\to0,\ \text{ as }k\to\infty.
\end{equation}
In order to prove \eqref{eq:conv:0}, we note that for any $\varepsilon>0$,
\begin{align*}
\E\int_0^s\left|{\bf1}_{H_r>t} - {\bf1}_{H^k_r>t}\right|dr&\le
\int_0^s\P(|H_r-H^k_r|>\varepsilon)dr+\E\int_0^s{\bf1}_{\{t-\varepsilon\le H^k_r\le t+\varepsilon\}}dr\\
&=\int_0^s\P(|H_r-H^k_r|>\varepsilon)dr+\E\int_{t-\varepsilon}^{t+\varepsilon}L^u_k(s)du.
\end{align*}
The first term on the right does not depend upon $t$ and tends to $0$ as $k\to\infty$ as a consequence of Corollary \ref{cor:estim}, while the second term is dominated by
\[ 2\varepsilon \ \sup_{t>0, \ k\ge1}L^t_k(s).\]
Hence \eqref{eq:conv:0} follows from Lemma \ref{le:estimLT} and the fact that $\varepsilon>0$ is arbitrary.
The Proposition is established.
\epf

We have in fact proved
\begin{corollary}\label{cor:convlt}
For any $t,s>0$, as $k\to\infty$,
\[ L^t_k(s)\to L^t(s)\ \text { in probabiity.}\]
\end{corollary}

We next establish
\begin{lemma}\label{contLs}
The local time $L^t(s)$ is continuous in $s$, for all $t\ge0$.
\end{lemma}
\bpf 
An argument very similar to that at the end of the proof of Proposition \ref{convH} yields the continuity of the map $s\mapsto\beta(H_s-t)^+-L^t(s)$, while the same Proposition implies that  $s\mapsto\beta(H_s-t)^+$ is continuous. The result follows.
\epf

Before proving the next Proposition, we show a uniform $L^p$-bound for the (approximating) local time(s), up to the time of the first big jump of $X$. To prepare this, we first fix $k\ge1$ and $s>0$, and consider the process
\[ A_t=\int_0^s\int_0^{\eps_k}{\bf1}_{H_r\le t} z\tilde{N}(dr,dz).\]
Let $\mathcal{G}_t$ denote the $\sigma$--algebra generated by the random variables
\begin{equation}\label{Ig}
 I_g=\int_0^s\int_0^{\eps_k}g(r,z)\, z\tilde{N}(dr,dz),
 \end{equation}
where $g$ is bounded and $\mathcal{P}\otimes\mathcal{B}_+$ measurable ($\mathcal{P}$ stands for the $\sigma$--algebra of predictable 
subsets of $\Omega\times\R_+$) and satisfies $\{g(r,z)=0\}\supset\{H_r>t\}$. We first establish
\begin{lemma}\label{le:mart}
The process $\{A_t: t\ge 0\}$ is a $(\mathcal{G}_t)$--martingale.
\end{lemma}
\bpf
It suffices to verify that $\E[(A_{t'}-A_t)I_g]=0$ for $t<t'$ and any $g$ as above, where $I_g$ is defined by \eqref{Ig}. This, however, 
is obvious. \epf

For $K>0$, let $\tau_K$ be the time of the first jump of $X$ of size greater than or equal to $K$.
We shall need the 
\begin{lemma}\label{le:boundLT} 
For any $p\ge1$, $s>0$, $K>0$, there exists a constant $C$ which depends only on those three parameters, such that 
\begin{align*}
 \sup_{t\ge0}\E[L^t(s\wedge \tau_K)^p]&\le C,\\
  \sup_{t\ge0, k\ge1}\E[L_k^t(s\wedge \tau_K)^p]&\le C.
  \end{align*}
\end{lemma}
\bpf
We shall prove the first inequality only. The second one is proved in exactly the same way.
Since $H_s$ and $L^t(s)$ are continuous in $s$, we can rewrite their expressions \eqref{betaH} and \eqref{defLT} as
\begin{align*}
\beta H_s&=X_{s-}-\inf_{0\le r<s} X_r-\int_0^{s-}\int_0^\infty(z+\inf_{r\le u< s}X_u-X_r)^+N(dz,dr),\\
L^t(s)&=\beta(H_s-t)^+\!\!-\!\!\int_0^{s-}{\bf1}_{H_r>t}dX_r
+\int_0^{s-}\int_0^\infty\!\!\!\!\!\!{\bf1}_{H_r>t}(z\!+\!\!\inf_{r\le u<  s}X_u-X_r)^+N(dr,dz).
\end{align*}
We first note that, since ${\bf1}_{[0,K]}(z)(z^2\wedge z^p)$ is $\pi$--integrable for all $p\ge1$,
one can easily show that for all $p\ge1, K>0, s>0$, there exists a constant $C_{p,K,s}$ such that
\begin{equation}\label{estimX}
\E\left(\sup_{0\le r<s\wedge \tau_K}|X_r|^p\right)\le C_{p,K,s}.
\end{equation}

We now estimate the last term in the above right hand side. It is clear that (the second inequality follows by combining $\beta H_s\ge0$ with the above identity)
\begin{align*}
\int_0^{s-}\int_0^\infty\!\!\!\!\!\!\!{\bf1}_{H_r>t}(z\!+\!\!\inf_{r\le u< s}X_u-X_r)^+N(dr,dz)&\le
\int_0^{s-}\int_0^\infty\!\!\!\!\!\!(z\!+\!\!\inf_{r\le u< s}X_u-X_r)^+N(dr,dz)\\
&\le X_{s-}-\inf_{0\le r< s}X_r\\
&\le 2\sup_{0\le r<s}|X_r|.
\end{align*}
 Next we observe that
\begin{align*}
\beta(H_s-t)^+&\le \beta H_s\le X_{s-}-\inf_{0\le r<s} X_r\\
&\le 2\sup_{0\le r<s}|X_r|.
\end{align*}
From the last two inequalities,
\begin{equation}\label{est1}
\begin{split}
 \beta(H_{s\wedge \tau_K}-t)^+ &+\int_0^{(s\wedge \tau_K)-}\!\int_0^\infty\!\!\!\!\!\!{\bf1}_{H_r>t}(z\!+\!\!\inf_{r\le u< s}X_u-X_r)^+N(dr,dz) \\
 &\le 4\sup_{0\le r<s\wedge\tau_K}|X_r|.
 \end{split}
\end{equation}
We now consider the second term
\begin{align*}
-\int_0^{s-} {\bf1}_{H_r>t} dX_r
=-\sqrt{2\beta}\int_0^s{\bf1}_{H_r>t} dB_r-\int_0^{s-}\int_0^\infty{\bf1}_{H_r>t}z\tilde{N}(dr,dz).
\end{align*}
The $p$--th absolute moment of the first term on the right is easy to estimate, since by the Burkholder--Davis--Gundy inequality,
\begin{equation}\label{est2}
\E\left(\left|\int_0^s{\bf1}_{H_r>t} dB_r\right|^p\right)\le C_p s^{p/2}.
\end{equation}
We finally estimate the $p$--th absolute moment of the last term.
Let $Y_s=\int_0^s\int_0^\infty{\bf1}_{H_r>t}z\tilde{N}(dr,dz)$. 
We first note that 
\begin{equation}\label{est3}
 |Y_{(s\wedge \tau_K)-}|\le\sup_{r\le s}\left|\int_0^r\int_0^K{\bf1}_{H_u>t}z\tilde{N}(du,dz)\right|.
 \end{equation}
Newt we use the Burkholder--Davis--Gundy inequality 
for possibly discontinuous martingales, see e.g. Theorem IV.48 in Protter \cite{Pr}, which yields
\begin{equation}\label{est4}
\begin{split}
\E\left(\sup_{r\le s}\left| \int_0^r\int_0^K {\bf1}_{H_u>t}z\tilde{N}(du,dz)\right|^p\right)&\le c_p
\E\left[\left(\int_0^{s}\int_0^K z^2N(du,dz)\right)^{p/2}\right].
\end{split}
\end{equation}
The result follows from a combination of \eqref{estimX}, \eqref{est1}, \eqref{est2}, \eqref{est3}, \eqref{est4}
and the fact that if $N$ is a Poisson random measure with mean measure $\nu$ and 
$f\in L^1(\nu)\cap L^\infty(\nu)$, then all moment of $N(f)$ are finite. The last statement can be deduced from the fact that 
the $k$--th cumulant of $N(f)$ is given as $\kappa_k(N(f))=\int f^k d\nu$, which is easy to verify for any step function $f$. 

 Just as the reflection of $H$ above zero leads to $L^0(s)>L^{0-}(s)=0$, the process $t\to L^t(s)$ is discontinuous, due to the fact that the jumps of $X$ create accumulations of local time of $H$ at certain level.
The points of discontinuity of $t\to L^t(s)$ are of course at most countable. They can be described
as follows : Let 
\[\mathcal{N}_s:=\{0\le r\le s;\ N(\{r\}\times \mathbb R_+) > 0\}\]
be the projection onto the $s$-axis of the support of the Poisson random measure $N$.
The set $\mathcal{N}_s$ is at most countable, and $\{H_r,\ r\in\mathcal{N}_s\}$ is the set of the points of discontinuity  of the mapping $t\to L^t(s)$.

\begin{proposition}\label{pro:contLT}
The local time $L^t(s)$ has  a version which is a.s. continuous in $s$ and cadlag in $t$.
\end{proposition}

\bpf

The continuity in $s$ has been established in Lemma \ref{contLs}.
Considering now the regularity in $t$, we note that the first term in the right of \eqref{defLT} is clearly continuous in $t$. Concerning the second term, we have for any $p>2$
and $t<t'$ from the Burkholder--Davis--Gundy inequality, with the stopping times introduced just before Lemma \ref{le:boundLT}, 
exploiting Jensen's inequality for the last inequality,\begin{align*}
\E&\left(\sup_{0\le r\le s\wedge \tau_K}\left|\int_0^r({\bf1}_{H_r>t}-{\bf1}_{H_r>t'})dB_r\right|^{2p}\right)\\
&\qquad\le4\E\left\{\left|\int_0^{s\wedge \tau_K}{\bf1}_{t<H_r\le t'}dr\right|^p\right\}\\
&\qquad=4(t'-t)^p\E\left\{\left(\int_t^{t'}L^u(s\wedge \tau_K)\frac{du}{t'-t}\right)^p\right\}\\
&\qquad\le 4(t'-t)^{p-1}\E\int_t^{t'}(L^u(s\wedge \tau_K))^pdu.
\end{align*}
This combined with Kolmogorov's lemma implies that the mapping $t\to \int_0^s{\bf1}_{H_r>t}dB_r$ has a version which is continuous in the two variables $t$ and $s$. 

Concerning the two last terms, if we replace the integrals over $(0,s]\times(0,\infty)$ by integrals over $(0,s]\times(\eps_k,\infty)$, then the sum of those two terms is c\`adl\`ag, the evolution between the jumps being absolutely continuous in the first term and decreasing in the second one. It remains to show that
the supremum over $t$ of
\[-\int_0^s\int_0^{\eps_k}{\bf1}_{H_r>t}z\tilde{N}(dr,dz)+\int_0^s\int_0^{\eps_k}\!\!\!\!\!{\bf1}_{H_r>t}(z+\inf_{r\le u\le s}X_u-X_r)^+N(dr,dz)\] tends to $0$ as $k\to\infty$. Concerning the second term, this follows from the fact that
\begin{align*}
\sup_t&\int_0^s\int_0^{\eps_k}\!\!\!\!\!{\bf1}_{H_r>t}(z+\inf_{r\le u\le s}X_u-X_r)^+N(dr,dz)\\
&\quad\le\int_0^s\int_0^{\eps_k}\!\!\!\!\!(z+\inf_{r\le u\le s}X_u-X_r)^+N(dr,dz),
\end{align*}
and the right hand side converges to 0 in probability as $k\to\infty$. 

Finally the uniform convergence of the first term follows from Lemma \ref{le:mart} and Doob's maximal inequality. \epf

\subsubsection{The Ray--Knight theorem}
We can now establish the Ray--Knight Theorem. 

\begin{theorem}\label{th:RNlin}
Under the assumption \eqref{recur} the stopping time $S_x$ defined in ~\eqref{defSx} is finite  a.s. and the process $\{L^t(S_x),\ t\ge0\}$ is
a CSBP with branching mechanism $\psi$.
\end{theorem}
\bpf
Let $S^k_x:=\inf\{s>0,\ L^0_k(s)>x\}$. Proposition \ref{RNcpp} shows that $\{L^t_{k}(S^k_x),\ t\ge0\}$ is
a CSBP with branching mechanism $\psi_k$ (here again $L_k$ denotes the local time of $H^k$).
It is plain that for any $g\in C_b(\R_+;\R_+)$ with compact support, we have both
\begin{align*}
\int_0^\infty g(t)L_k^t(S^k_x)dt&=\int_0^{S_x^k}g(H^k_s)ds,\quad\text{and}\\
\int_0^\infty g(t)L^t(S_x)dt&=\int_0^{S_x}g(H_s)ds.
\end{align*}
 Provided we show that $S^k_x\to S_x$, which will be done in the next Lemma, it follows from Corollary \ref{cor:convH} that the right--hand side of the first identity converges in probability to the right--hand side  of the second identity
 in probability, as $k\to\infty$. Consequently   for any $T>0$,
 \[ L_k^\cdot(S_x)\to L^\cdot(S_x)\]
in $L^2(0,T)$ weakly, in probability, as $k\to\infty$.

On the other hand, from Proposition \ref{RNcpp}, $\{L_k^t(S^k_x),\ t\ge0\}$ is a CSBP with branching mechanism $\psi_k$. Let now $W$ be a space--time white noise, and $M$ a Poisson random measure with mean $ds\times du\times\pi(dz)$,
while $\widetilde{M}$ will denote the compensated measure $M(ds,du,dz)-ds\,du\,\pi(dz)$.
It is clear that if $\{Z^{k,x}_t,\ t\ge0\}$ denotes the unique strong solution of the Dawson--Li type SDE (see \cite{DL})
\begin{equation*}
\begin{split}
Z^{k,x}_t&=x+\alpha\int_0^t Z^{k,x}_sds+\sqrt{2\beta}\int_0^t\int_0^{Z^{k,x}_s}W(ds,du)\\
&\quad+\int_0^t\int_0^{Z^{k,x}_{s-}}\int_{\eps_k}^\infty z\widetilde{M}(ds,du,dz),
\end{split}
\end{equation*}
then for each $k\ge1$, $\{L_k^t(S^k_x),\ t\ge0,\ x>0\}$ and $\{Z^{k,x}_t,\ t\ge0,\ x>0\}$
have the same law. On the other hand, it is not hard to show that $Z^{k,x}_t\to Z^x_t$
in probability, locally uniformly in $t$, where $Z^x_t$ is the unique solution of the SDE
\begin{equation}\label{eq:Z-CSBP}
\begin{split}
Z^{x}_t&=x+\alpha\int_0^t Z^{x}_sds+\sqrt{2\beta}\int_0^t\int_0^{Z^{x}_s}W(ds,du)\\
&\quad+\int_0^t\int_0^{Z^{x}_{s-}}\int_{0}^\infty z\widetilde{M}(ds,du,dz).
\end{split}
\end{equation}
The result follows from a combination of the above arguments. \epf

It remains to show that $S^k_x\to S_x$.
\begin{lemma}\label{le:convS}
For any $x>0$, as $k\to\infty$,
\[ S^k_x\to S_x\ \text{ in probability}.\]
\end{lemma}
\bpf
From the definition of $S_x:=\inf\{s>0,\ L^0(s,H)>x\}$, for any $\eps>0$, $L^0(S_x+\eps)>x$.
Hence $\limsup_{k\to\infty}S^k_x\le S_x$. However, $L^0(s,H)=-\inf_{0\le r\le s}X_r$. By Theorem VII.1 of
Bertoin \cite{B}, the process $x\mapsto S_x$ is a subordinator. Consequently, by Proposition I.7 of \cite{B}, a.s. 
$S_x=S_{x-}=\inf\{s>0,\ L^0(s,H)\ge x\}$. So for any $\eps>0$, $L^0(S_x-\eps,H)<x$, and a.s. $\liminf_{k\to\infty}S^k_x\ge S_x$.
\epf

\section{The case with interaction}\label{sec:interact}
For the rest of the paper we consider, instead of \eqref{eq:CSBP}, the collection of SDE's~\eqref{eq:Zinteract}. In other words, the  linear drift term  $-\alpha Z_t^x \, dt$ in \eqref{eq:CSBP}  is replaced by the non-linear drift term $f(Z_t^x) dt$,  with $f$ satisfying \eqref{hyp:f}.

Connecting to the results of the previous section, we consider  a process $X$ defined by  \eqref{X} with $\alpha=0$, i.e.
%
\begin{equation}\label{Xsimplified}
X_s=\sqrt{2\beta} B_s+\int_0^s\int_0^\infty z\tilde{N}(dr,dz), \quad s\ge 0, 
\end{equation}
where again $\tilde{N}$ denotes the compensated measure $\tilde{N}(dr,dz)=N(dr,dz)-dr\, \pi(dz)$.

Our final aim in this paper is to obtain a Ray--Knight representation for the solution $Z$ of  \eqref{eq:Zinteract} in terms of an appropriate height process. 
For this, our strategy will be to introduce, via Girsanov's theorem, the appropriate drift into the equation  \eqref{betaH} for the height process $H$. This change of measure will introduce the same drift into the process $X$, and should lead to the SDE's \eqref{eq:H} and \eqref{eq:X} for the pair $(X,H)$. 

However, condition \eqref{hyp:f} guarantees only {\em local} boundedness of $f'$. Thus, in order to make sure that Girsanov's theorem is applicable, we use a localization procedure and associate to each $b \in (0,\infty)$ a function
\begin{equation}\label{def:fb}
\begin{split}
 f_b&\in C^1_b(\R_+),\  f'_b \text{ is uniformly continuous on }\R_+,\\
 \text{and } f_b(z)&=f(z),\ 0\le z\le b\,.
 \end{split}
 \end{equation}
We also assume that $f'_b(z)\le\theta$, for all $z>0$, $b>0$.

Even with this localization,  the process $H$ (which then solves \eqref{eq:H} with $f_b$ instead of $f$) might tend to infinity before its local time at $t=0$ has achieved the value $S_x$, $x >0$. Then there would be no way to make sense of the process $L^t(S_x)$. One way to circumvent this difficulty would be
 to define $H$ reflected below an arbitrary level~$a$ as in \cite{JFD} and \cite{P}, and identifying the law of $L^\cdot(S_x)$ as that
 of $Z^x$, killed at time $t=a$. However, there would be difficulties with the definition of the thus reflected SDE for $H$, due to the jump terms. Therefore, we will use an additional localization by adding a drift which acts only while $H$ takes values above $a > 0$, and has the effect of forcing $H$ to hit $0$ after any time $s_0>0$, i.e. $\inf\{s>s_0, H_s=0\}<\infty$ a.s.. Our choice for this will be
 \[ g_a(h)=-(h-a)^+.\]
 After taking the limit $b\to \infty$ we will identify the law of $\{L^t(S_x),\, 0\le t\le a\}$ with that of $\{Z^x_t,\, 0\le t\le a\}$, but will loose the interpretation of  $L^t(S_x)$ for $t>a$.
\subsection{The  case $\pi=0$}\label{intpi0}
This case is treated in Pardoux \cite{P}.  The equation \eqref{eq:H} for the height process $H$ reads
\begin{equation}\label{SDEHint}x
\beta H_s=\int_0^sf'\left(L^{H_r}(r)\right)dr+\sqrt{2\beta}B_s-\inf_{0\le r\le s}X_r.
\end{equation}
It is shown in Proposition 16 of \cite{P} that in this case the process $\{Z^x_t,\ t\ge0\}$ goes extinct
a.s. in finite time for all $x>0$ iff
\begin{equation}\label{condit:subcrit}
  \int_1^\infty \exp\left(-\frac{1}{\beta}\int_1^u \frac{f(r)}{r}dr\right)du=\infty,
  \end{equation}
and in this case Corollary 7 of \cite{P} shows that $\{L_{S_x}^t,\ t\ge0\}$ solves the
SDE \eqref{eq:Zinteract}, of course with $\widetilde{M}\equiv0$.  If the condition \eqref{condit:subcrit}
is not satisfied, $Z_t^x$ need not go extinct, and in that case the process $H_s$ may tend to infinity as $s\to\infty$,
so that we may have $L^0(\infty)<x$. However one can still obtain an extension of the second Ray--Knight theorem, by reflecting $H$ below an arbitrary level, as in Delmas \cite{JFD}, see Theorem 14
in \cite{P}.  The equation \eqref{SDEHint} has a unique weak solution: for each $x > 0$, existence up to time $S_x$   follows from Girsanov's theorem, see the explanation on p.~95-97 together with Corollary 8 in \cite{P}. Since Girsanov's theorem can be applied also in the reverse direction, this implies weak uniqueness up to time $S_x$; see also \cite{PW} Sec.~4.1 for that argument in the case of affine linear $f'$.

\subsection{The  case of finite $\pi$}
In this subsection, we assume that $\pi((0,+\infty))<\infty$.
We now use Girsanov's theorem in order to describe the corresponding height process. Under the reference measure $\P$, let $H$ denote the solution of \eqref{betaH}.
For any  $a, b> 0$, let $Y^{a,b}$ denote the following Girsanov Radon--Nikodym derivative 
\[Y^{a,b}_s=\exp\left(\frac{1}{\sqrt{2\beta}}\int_0^{s\wedge S_{x}}[f_b'(L^{H_r}(r))+g_a(H_r)]dB_r-\frac{1}{4\beta}\int_0^{s\wedge S_{x}} [f_b'(L^{H_r}(r))+g_a(H_r)]^2dr\right),\]
and define
 \begin{equation}\label{defBMa}
  B^{a,b}_s=B_s-\frac{1}{\sqrt{2\beta}}\int_0^s[f_b'(L^{H_r}(r))+g_a(H_r)]dr, \quad s \ge 0.
  \end{equation}
This is a Brownian motion up to time $S_x$ under the unique probability measure $\P^{a,b}$ which is such that, with $\F_s=\sigma\{H_r,\ 0\le r\le s\}$,
 \begin{equation}\label{RNderiv}
  \frac{d\P^{a,b}}{d\P}\Big|_{\F_s}=Y^{a,b}_s,\ s>0.
  \end{equation}
 Since $f'_b$ and $g_a$ are bounded, this follows readily from Proposition 35 in \cite{P}. It is easy to verify that the law of the random measure 
 $N(dr,dz)$ is the same under $\P^{a,b}$ and under $\P$. Indeed, one way to check that under $\P^{a,b}$, $N$ is a Poisson random measure on $(0,+\infty)^2$ with mean the Lebesgue measure is to check that for any $\varphi\in C((0,+\infty)^2;\R_+)$ with compact support, 
\begin{equation}\label{PPP}
\E^{a,b}\exp[N(-\varphi)]=\exp\{\int_0^\infty\int_0^\infty[e^{-\varphi(s,z)}-1]dsdz\},
\end{equation}
where we have used the notation
 \[ N(\psi)=\int_{(0,\infty)^2}\psi(r,z)N(dr,dz)\,.\]
To verify \eqref{PPP}, note first that
 it follows from It\^o's formula that for any $s>0$,
 if $\varphi_s(r,z):=\varphi(r,z){\bf1}_{[0,s]}(r)$,
  \begin{align*}
 \exp[N(-\varphi_s)]&=1+\int_0^s \exp[N(-\varphi_{r-})]\int_0^\infty\left(e^{-\varphi(r,z)}-1\right)\tilde{N}(dr,dz)\\ &\qquad
 +\int_0^s \exp[N(-\varphi_{r})]\int_0^\infty\left(e^{-\varphi(r,z)}-1\right)dzdr\, .
  \end{align*}
  We now observe that  the above integral with respect to $\tilde{N}$ is a martingale under $\P^{a,b}$, which follows,
see e.g. Theorem III.36 in Protter \cite{Pr},  from the fact that both it is a martingale under $\P$, and its quadratic covariation with the Radon-Nikodym derivative \eqref{RNderiv} vanishes, i.e. 
 \[ \langle Y^{a,b}_\cdot,\int_0^\cdot \exp[N(-\varphi_{r-})]\int_0^\infty\left(e^{-\varphi(r,z)}-1\right)\tilde{N}(dr,dz)\rangle\equiv0.\]
 This readily implies that
  \begin{align*}
 \E^{a,b} \exp[N(-\varphi_s)]=1+\int_0^s \E^{a,b}\exp[N(-\varphi_{r})]\int_0^\infty\left(e^{-\varphi(r,z)}-1\right)dzdr,
 \end{align*}
 from which \eqref{PPP} follows by explicit integration of a linear ODE, choosing $s$ large enough so that supp$(\varphi)\subset[0,s]\times(0,+\infty)$. 
  
 It follows from \eqref{Xsimplified} and \eqref{defBMa} that
 \begin{align*} 
 X_s&=\int_0^s[f_b'(L^{H_r}(r))+g_a(H_r)]dr+X^{a,b}_s,\quad  s\ge0.
\end{align*}
where
 \begin{align*} 
 X^{a,b}_s=\sqrt{2\beta}B^{a,b}_s+\int_0^s\int_0^\infty z\tilde{N}(dr,dz).
\end{align*}
Consequently \eqref{betaH} can be written as
\begin{equation}\label{SDEH}
\beta H_s=\int_0^s[f_b'(L^{H_r}(r))+g_a(H_r)]dr+X^{a,b}_s-\inf_{0\le r\le s}X_r-\int_0^s\int_0^\infty\left(z+\inf_{r\le u\le s}X_u-X_r\right)^+N(dr,dz).
\end{equation}

 Weak existence of a solution  to \eqref{SDEH} follows from the above explicit 
 construction. Weak uniqueness follows from the fact that \[\frac{d\P}{d\P^{a,b}}\Big|_{\F_s}=(Y^{a,b}_s)^{-1}.\]
We denote again by $L^t(s)$ the local time accumulated by the process $H$ at level $t$ up to time $s$, and $S_x=\inf\{s>0;\, L^0(s)>x\}$.

We have
\begin{lemma}
For any $a, b > 0$ we have $\P^{a,b}(S_x<\infty)=1$.
\end{lemma}
\bpf 
We observe that 
\[f_b'(L^{H_r}(r))+g_a(H_r)\le\theta-(H_r-a)^+.\] 
Consequently, whenever $H_r>\theta+a+1$, the drift in the equation for $H$ is bounded from above by $-1$. This is enough to conclude that under $\P^{a,b}$, the process $H$ returns to $0$ after arbitrarily large times, hence accumulates arbitrary quantities of local time at level $0$.
\epf

We can rewrite \eqref{SDEH} as
 \begin{align}\label{SDEHfpi}
 \begin{split}
 \beta H_s&=\int_0^s\left[f_b'(L^{H_r}(r))+g_a(H_r)-\gamma\right]dr+\sqrt{2\beta}B^{a,b}_s-\inf_{0\le r\le s}X_r \\
 &\quad+ \int_0^s\int_0^\infty\left(z-\left[z+\inf_{r\le u\le s}X_u-X_r\right]^+\right)N(dr,dz),
 \end{split}
 \end{align}
 where $\gamma=\int_{(0,+\infty)}z\, \pi(dz)$.
 \begin{proposition}\label{RNcpp_f}
Assume that the measure $\pi$ is finite, and fix $a,b >0$. Under $\P^{a,b}$, the process $\{L^t(S_{x}),\ 0\le t\le a,\ x>0\}$ is, on the time interval $[0,a]$,
a solution of the collection indexed by $x>0$ of SDEs
\begin{equation}\label{eq:Zinteractb}
\begin{split}
Z^{x,b}_t&=x+\int_0^t f_b(Z^{x,b}_r)dr+\sqrt{2\beta}\int_0^t\int_0^{Z^{x,b}_r}W(dr,du)\\
&\qquad+\int_0^t\int_0^{Z^{x,b}_{r-}}\int_{0}^\infty z\widetilde{M}(dr,du,dz), \quad t\ge0.
\end{split}
\end{equation}
\end{proposition}
\bpf
{\sc Step 1. Equation for $L^t(S_x)$} Here $x>0$ is fixed. We first note that $H_{S_x}=0$ implies that 
$\sum_{0\le r\le S_x}\Delta \overline{X}^{S_x}_r=0$. Moreover $X_{S_x}=-x$. Consequently formula \eqref{defLT} at $s=S_x$ 
reads
\begin{align*}
L^t(S_x)&= -\int_0^{S_x}{\bf1}_{H_r>t}dX_r\\
&=x+\int_0^{S_x}{\bf1}_{H_r\le t}dX_r\, .
\end{align*}
Note that
\begin{align*}
X_s&=B_s+\int_0^s\int_0^\infty z\tilde{N}(dr,dz)\\
&=\int_0^s[f_b'(L^{H_r}(r)+g_a(H_r)]dr+B^{a,b}_s+\int_0^s\int_0^\infty z\tilde{N}(dr,dz),\ \text{ hence for $t\le a$}\\
L^t(S_x)&=x+\int_0^{S_x}{\bf1}_{H_r\le t}f_b'(L^{H_r}(r))dr +\int_0^{S_x}{\bf1}_{H_r\le t}dB^{a,b}_r 
+\int_0^{S_x}{\bf1}_{H_r\le t}\int_0^\infty z\tilde{N}(dr,dz)\, ,
\end{align*}
where we have exploited the fact that for $t\le a$, ${\bf1}_{H_r\le t}g_a(H_r)\equiv0$.
We will now rewrite each of the three integrals of the last right hand side. For the first one, we use, similar to an argument on page 728 of \cite{PW}, the generalized occupation times formula from Exercise 1.15 in Chapter VI of \cite{RY}, and obtain
\begin{align*}
\int_0^{S_x}{\bf1}_{H_r\le t}f_b'(L^{H_r}(r))dr&=\int_0^t\int_0^{S_x}f_b'(L^u(r))dL^u_r du\\
&=\int_0^tf_b(L^u(S_x))du\,.
\end{align*}
We next consider the process
\[ U_t:=\int_0^{S_x}{\bf1}_{H_r\le t}\,dB^{a,b}_r, \quad t\ge 0.\]
For $t\ge 0$, let $\HH^B_t$ denote the sigma--algebra generated by  the random variables of the form
\[ Y_g=\int_0^{S_x}g(r)dB^{a,b}_r,\]
where $g$ is progressively measurable and satisfies $\{g(r)=0\}\supset\{H_r>t\}$. It is easily seen that $U=(U_t)_{t\ge 0}$  is an $\HH^B$--martingale for the filtration $\HH^B= (\HH^B_t)_{t\ge 0}$. We now show that it is a continuous martingale. Indeed, for any $K>0$, let again $\tau_K$ denote the time of the first jumps of $X$ is size greater than $K$. On the event $\Omega_{x,K}=\{S_x\le\tau_K\}$, for any $t>0$, 
\[ U_t=\int_0^{S_x\wedge\tau_K}{\bf1}_{H_r\le t}dB^{a,b}_s.\]
Therefore, for $t'>t>0$, $p>2$,
\begin{align*}
\E\left[|U_{t'}-U_{t}|^p;\Omega_{x,K}\right]&=\E\left[\left| \int_0^{S_x\wedge\tau_K}{\bf1}_{t<H_r\le t'}dB^{a,b}_r\right|^p\right]\\
&=\E\left[\left| \int_0^{S_x\wedge\tau_K}{\bf1}_{t<H_r\le t'}dr\right|^{p/2}\right]\\
&=\E\left[\left| \int_t^{t'}L^u(S_x\wedge\tau_K)du\right|^{p/2}\right]\\
&\le \sup_{u>0}\E\left(\left|L^u(S_x\wedge\tau_K)\right|^{p/2}\right)\times|t'-t|^{p/2}.
\end{align*}
The a.s. continuity of $U$ follows from this computation, Lemma \ref{le:boundLT}, Kolmogorov's Lemma, and the fact that 
$\P\left(\cup_{K\ge1}\Omega_{x,K}\right)=1$.

We next note that
\[ \langle U\rangle_t=\int_0^tL^u(S_x)du \,, \quad t\ge 0.\]
Indeed, by  It\^o's formula,
\[ U_t^2-\int_0^tL^u(S_x)du=2\int_0^{S_x} {\bf1}_{H_r\le t}\int_0^r{\bf1}_{H_s\le t}dB^{a,b}_s dB^{a,b}_r, \quad t\ge 0,\]
is a $\HH^B$--martingale.

It is now clear that there exists a space--time white noise $W(dr,du)$ such that
\[ U_t=\int_0^t\int_0^{L^r(S_x)}W(dr,du), \quad t\ge 0.\]
Here the choice of representing the above martingale as a stochastic integral with respect to a space time white noise, rather than with respect to a Brownian motion, is motivated by Step 2 of the proof below.

We finally consider the process
\begin{equation*}\label{decompV}
 \int_0^{S_x}{\bf1}_{H_r\le t}\int_0^\infty z\tilde{N}(dr,dz)=V_t- \gamma\int_0^t L^u(S_x)du, \quad t\ge 0,
\end{equation*}
where
\[ \gamma=\int_0^\infty z\pi(dz),\quad\text{and }\quad V_t:=\int_0^{S_x}{\bf1}_{H_r\le t}\int_0^\infty z{N}(dr,dz)\, ,\]
so that we have obtained
\begin{equation}\label{eq:LSx}
L^t(S_x)=x+\int_0^t[f_b(L^r(S_x))-\gamma L^r(S_x)]dr+\int_0^t\int_0^{L^r(S_x)}W(dr,du)+V_t\,.
\end{equation}
We now want to rewrite the process $V_t$ in a different way. For that sake, we use again the construction introduced in Proposition
\ref{RNcpp}. 

We start with $X^{(0)}, H^{(0)}, L_{(0)},S^{(0)}_x$ defined as follows.
\begin{align*}
X^{(0)}_s&=\int_0^s[f_b'(L^{H^{(0)}_r}_{(0)}(r))-\gamma]dr+\sqrt{2\beta} B^{a,b}_s,\\
\beta H^{(0)}_s&=\int_0^s[f_b'(L^{H^{(0)}_r}_{(0)}(r))-\gamma]dr+\sqrt{2\beta} B^{a,b}_s-\inf_{0\le r\le s}X^{(0)}_r\, ,\\
L^t_{(0)}(s)\ &\text{ is the local time accumulated by }H^{(0)}\ \text{ at level }t\ \text{ up to time }s\, ,\\
S^{(0)}&=\inf\{s>0,\ L^0_{(0)}(s)>x\}\, .
\end{align*}
Let $N^{(0)}$ denote an independent copy of the Poisson random measure $N$, and \mbox{$\{(s_i,z_i),\ 1\le i\le J\}$} be the set of points of  $N^{(0)}$ on $[0,S^{(0)}_x]\times(0,+\infty)$. If $J=0$, then 
$(X,H,L,S_x)\equiv(X^{(0)},H^{(0)},L_{(0)},S^{(0)}_x)$, and we are done. Otherwise, we select the a.s.~unique index $m \in \{1,\ldots, J\}$ such that $H^{(0)}_{s_m}=\min_{1\le i\le J}H^{(0)}_{s_i}$, and we define $X^{(1)}, H^{(1)}, L_{(1)},S^{(1)}_x$ as follows. We consider an independent copy ${B}^{a,1}$ of $B^{a}$, and define
\begin{align*}
X^{(1)}_s&=\begin{cases} X^{(0)}_s,&\text{for $s\le s_m$},\\
X^{(0)}_{s_m}+z_m+\int_{s_m}^s[f_b'(L^{H^{(1)}_r}_{(1)}(r))-\gamma]dr+\sqrt{2\beta} B^{a,1}_{s-s_m},&\text{for $s_m<s\le s_m+\tilde{s}_1$},\\
X^{(0)}_{s_m}+\int_{s_m+\tilde{s}_1}^s[f_b'(L^{H^{(1)}_r}_{(1)}(r))-\gamma]dr+\sqrt{2\beta} B^{a,b}_{s-\tilde{s}_1},&\text{for $s\ge s_m+\tilde{s}_1$}\end{cases}\\
\beta H^{(1)}_s&=\int_0^s[f_b'(L^{H^{(1)}_r}_{(1)}(r))-\gamma]dr+\sqrt{2\beta} B^{a,b}_s-\inf_{0\le r\le s}X^{(1)}_r\, ,\\
L^t_{(1)}(s)\ &\text{ is the local time accumulated by }H^{(1)}\ \text{ at level }t\ \text{ up to time }s\, ,\\
S^{(1)}_x&=\inf\{s>0,\ L^0_{(1)}(s)>x\}\, ,
\end{align*}
where
\[ \tilde{s}_1=\inf\left\{s>0, X^{(1)}_{s_m+s}<X^{(0)}_{s_m}\right\}\, .\]
Note that $S^{(1)}_x=S^{(0)}_0+\tilde{s}_1$, since it is true under the reference probability $\P$ (see the same construction in Proposition \ref{RNcpp}).
We next define as follows the Poisson random measure $N^{(1)}$ on $[0,S^{(0)}_x]\times(0,+\infty)$. Given $\tilde{N}^{(1)}$  an independent copy of $N$, which we restrict to $[0,\tilde{s}_1]\times(0,+\infty)$, the points of $N^{(1)}$ are those of $N^{(0)}$ on
$[0,s_m]\times(0,+\infty)$, those of $\tilde{N}^{(1)}$ whose first coordinate has been shifted by $+s_m$ on 
$[s_m,s_m+\tilde{s}_1]\times(0,+\infty)$, and finally those of the restriction of $N$ to $[s_m,S^{(0)}_x]\times(0,+\infty)$ shifted by $+\tilde{s}_1$ on $[s_m+\tilde{s}_1,S^{(1)}_x]\times(0,+\infty)$. 

We are now ready to iterate our procedure, and construct the elements indexed by $2$. The iteration terminates a.s. at rank $K\ge J$ which is such that $N^{(K)}$ has no point. The law of $K$ is that of the number of points of our original Poisson random measure $N$ on $[0,S_x]\times(0,+\infty)$. Note that starting from $X,H,L,S_x$, we could construct a copy of the above sequence in reverse order by deleting one by one the jumps of $X$ on $[0,S_x]$, starting from the one corresponding to the largest value of 
$H$.

Coming back to the above sequence, the jumps of $\{V_t,\ t>0\}$ are described by that sequence in the order in which they appear as $t$ increases. It follows from our construction that the process $V$ can be written as
\[ V_t = \int_0^t\int_0^{L^{r-}(S_x)}\int_0^\infty zM(dr,du,dz), \quad t\ge 0,\]
where $M$ is a Poisson random measure on $(0,+\infty)^3$ with mean measure $dr\, du\, \mu(dz)$ as in Proposition \ref{RNcpp}.

Inserting this formula for $V$ in \eqref{eq:LSx}, we have proved that for fixed $x>0$, the process $\{L^t(S_x),\ 0\le t\le a\}$ satisfies equation \eqref{eq:Zinteract}.

\noindent{\sc Step 2 Identification of the law of $\{L^t(S_{x}),\, 0\le t\le a,x>0\}$}
If we define $\tilde{H}^x_s=H_{S_x+s}$ and $\tilde{X}^x_s=X_{S_x+s}+x$, we have $\tilde{X}^x_0=0$ and under $\P$
\[ \tilde{H}_s^x= \tilde{X}^x_s-\inf_{0\le r\le s}\tilde{X}^x_r+\int_0^s\int_0^\infty(z+\inf_{r\le u\le s}\tilde{X}^x_u-\tilde{X}^x_r)^+N(dr,dz)\,.\]
Denote again $\F_s=\sigma\{X_r,\,  0\le r \le s\}$.
It is not hard to see that under $\P^{a,b}$ $\{L^t(S_x+s)-L^t(S_x),\, s\ge0, 0\le t\le a\}$ is a function of both $\{L^{t'}(S_x),\,0\le t'\le a\}$ (through the nonlinear coefficient 
$f'_b$), and noises which are independent of $\FF_{S_x}$. Now we fix both $x$ and $y>0$, and note that
\begin{align*}
 L^t(S_{x+y})-L^t(S_x)&=y+\int_0^{S_{x+y}-S_x}{\bf1}_{\tilde{H}^x_r\le t}f_b'(L^{\tilde{H}^x_r}(S_x+r))dr\\&\quad
+\int_0^{S_{x+y}-S_x}\!\!\!\!\!\!\!\!{\bf1}_{\tilde{H}^x_r\le t}dB^{a,b}_{S_x+r}+\int_0^{S_{x+y}-S_x}\!\!\!\!\!\!\!\!{\bf1}_{\tilde{H}^x_r\le t}\int_0^\infty \!\!\!\!z\tilde{N}(S_x+dr,dz)\,.
\end{align*}
Applying the same extended occupation times formula as above, we deduce that 
\[ \int_0^{S_{x+y}-S_x}{\bf1}_{\tilde{H}^x_r\le t}f_b'(L^{\tilde{H}^x_r}(S_x+r))dr
= \int_0^t\left[f_b(L^u(S_{x+y}))-f_b(L^u(S_x)))\right]du\,. \]
From the same arguments of the previous steps, we see that $Z^{x,y,b}_t:= L^t(S_{x+y})-L^t(S_x)$ satisfies for $0\le t\le a$
\begin{align*}
Z^{x,y,b}_t&=y+\int_0^t\left[f_b(L_r(S_x)+Z^{x,y,b}_r)-f_b(L_r(S_x))\right]dr+
\sqrt{2\beta}\int_0^t\int_0^{Z^{x,y,b}_r}W^x(dr,du)\\
&\quad+\int_0^t\int_0^{Z^{x,y,b}_{r-}}\int_0^\infty z\widetilde{M}^x(dr,du,dz),
\end{align*}
where $W^x$ and $W$ (resp. $M^x$, $M$) are i.i.d. The independence follows by noting that the cross quadratic variation is zero. 
It follows from the independence property of the white noise and the Poisson random measure on disjoint subsets that the pair 
$\{(L^t(S_x),L^t(S_{x+y})-L^t(S_x)),\, 0\le t\le a\}$ has the same law as $\{(Z^{x,b}_t,Z^{x+y,b}_t-Z^{x,b}_t),\, 0\le t\le a\}$, hence also the two pairs $\{(L^t(S_x),L^t(S_{x+y})),\, 0\le t\le a\}$ and $\{(Z^{x,b}_t,Z^{x+y,b}_t),\, 0\le t\le a\}$ have the same law. 

A similar argument shows that for any $n\ge2$ and  $x_1<x_2<\cdots<x_n$, the two $n$--dimensional processes
$\{(L^t(S_{x_1}),L^t(S_{x_2}),\ldots,L^t(S_{x_n})),\, 0\le t\le a\}$
and $\{(Z^{x_1,b}_t,Z^{x_2,b}_t,\ldots,Z^{x_n,b}_t),\, 0\le t\le a\}$ have the same law. This proves the result.
\epf

\subsection{The  general case}

 With $\eps_k$ and   $\pi_k$ as in Section \ref{Sec241}, we are now going to take the limit as $k \to \infty$ in the setting of
the previous subsection. To this end, we first fix $a,b >0$. Since the drift $f_b'(L^{H_r}(r))$ is  not Lipschitz in $H$ with respect to any
of the standard  metrics on the continuous paths, 
it seems that the only practical route to access $H$ and its local time, and to establish our final result Theorem~\ref{RNfinal},  is to rely on the convergence
result of Section \ref{sec:gen_case}, and Girsanov's theorem. We recall that $|f_b'|$ and $g_a$ are bounded.

Consider the sequence $H^{k}$,  $k \ge 1$, of Section  \ref{sec:gen_case},  let $L_{k}$ denote the
 local time of $H^k$ an define $S^k_x=\inf\{s>0;\, L^0_k(s)>x\}$. We need to take the limit in the sequence of Radon--Nikodym derivatives
\[ Y^{a,b,k}_s=\exp\left(\frac{1}{\sqrt{2\beta}}\int_0^{s\wedge S^k_{x}}[f_b'(L_{k}^{H^{k}_r}(r))+g_a(H^k_r)]dB_r-
\frac{1}{4\beta}\int_0^{s\wedge S^k_x} |f_b'(L_{k}^{H^{k}_r}(r))+g_a(H^k_r)|^2dr\right).\]

The reference probability $\P$ governs the case $f_b\equiv0$ and $g_a\equiv 0$. Hence under  $\P$,
\begin{equation}\label{eq:Hak}
\beta H^{k}_s=X^k_s-\inf_{0\le r\le s}X^k_r-\int_0^s\int_{\eps_k}^\infty\left(z+\inf_{r\le u\le s}X^k_u-X^k_r\right)^+N(dz,dr).
\end{equation}
The quantities introduced in the previous subsection need now to be indexed by $k\ge1$. That is we consider the probability measure $\P^{a,b,k}$ such that for all $s>0$, with $\F^{k}_s=\sigma\{H^{k}_r,\ 0\le r\le s\}$,
\[ \frac{d\P^{a,b,k}}{d\P}\Big|_{\F^{k}_s}=Y^{a,b,k}_s,\ s>0.\]
Under $\P^{a,b,k}$, the process $H^{k}$ solves the SDE (see \eqref{SDEH})
\begin{equation}\label{eq:Hakf}
\begin{split}
\beta H^{k}_s&=\int_0^s \left[f_b'( L^{H^{k}_r}(r))+g_a(H^{k}_r)\right]dr+X^{a,b,k}_s-\inf_{0\le r\le s}X^k_r
 \\
 &\quad-\int_0^s\int_{\eps_k}^\infty\left(z+\inf_{r\le u\le s}X^k_u-X^k_r \right)^+N(dr,dz),
 \end{split}
 \end{equation}
 where \[X^{a,b,k}_s=\sqrt{2\beta} B^{a,b}_s+\int_0^s\int_{\eps_k}^\infty z\tilde{N}(dr,dz),\]
 and 
 \[ B^{a,b}_s=B_s-\frac{1}{\sqrt{2\beta}}\int_0^s[f_b'(L^{H^{k}_r}(r))+g_a(H^{k}_r)]dr\]
 is a Brownian motion under $\P^{a,b,k}$, up to time $S^k_x$.

  Under the reference probability $\P$, $H$ is defined by \eqref{betaH}, that is
 \begin{equation}\label{eqHa}
 \beta H_s=X_s-\inf_{0\le r\le s}X_r-\int_0^s\int_0^\infty\left(z+\inf_{r\le u\le s}X_u-X_r\right)^+N(dz,dr).
 \end{equation}
  The definition of the pair $(Y^{a,b},\P^{a,b})$, which was given at the beginning of the previous subsection for the case of a finite $\pi$, remains the same also for a general $\pi$ satisfying \eqref {hyp:pi}.
 Under $\P^{a,b}$, $H$ solves the SDE
 \begin{equation}\label{eq:Haf}
\begin{split}
\beta H_s&=\int_0^s \left[f_b'\left( L^{H_r}(r)\right)+g_a(H_r)\right]dr+X^{a,b}_s-\inf_{0\le r\le s}X_r\\
 &\quad-\int_0^s\int_0^\infty\left(z+\inf_{r\le u\le s}X_u-X_r \right)^+N(dr,dz).
 \end{split}
 \end{equation}
Again by the argument developed in the previous subsection, \eqref{eq:Haf} has a unique weak solution. The main argument in this subsection is
\begin{proposition}\label{limitGirsanov}
Let $a, b>0$ and $s> 0$ be fixed. Then, under the reference probability measure $\mathbb P$,
 $Y^{a,b,k}_s\to Y^{a,b}_s$ as $k\to\infty$, in probability and also in $L^p$ for any $p\ge 1$.
\end{proposition}
\bpf
Since $|f_b'|$ and $g_a$ are bounded, for any $p\ge1$, $\{(Y^{a,b,k}_s)^p\}_{k\ge1}$ is uniformly integrable, hence
it suffices to establish the convergence in probability. 
For that purpose, we need to show that
\begin{align*}
\int_0^{s} \left|{\bf1}_{r\le S_x}f_b'(L^{H_r}(r))-{\bf1}_{r\le S^k_x}f_b'(L_k^{H^{k}_r}(r))\right|^2dr&\to0\\
 \text{and }\quad \int_0^s\left|{\bf1}_{r\le S_x}g_a(H_r)-{\bf1}_{r\le S^k_x}g_a(H^k_r)\right|^2 dr&\to0,
 \end{align*}
 as $k\to\infty$ in $\mathbb P$-probability. 

The second convergence follows readily from Lemma \ref{le:convS}, Corollary \ref{cor:convH},  the Lipschitz continuity of $g_a$ and the dominated convergence theorem. The rest of this proof will be devoted to establishing the first convergence.

For this purpose we consider
\begin{align*} 
\left|{\bf1}_{r\le S_x}f_b'\left(L^{H_r}(r)\right)-{\bf1}_{r\le S^k_x}f_b'\left(L_k^{H^{k}_r}(r)\right)\right|&\le 
\left|f_b'\left(L^{H_r}(r)\right)-f_b'\left(L^{H^{k}_r}(r)\right)\right| \\
&\qquad+\left|{\bf1}_{r\le S_x}f_b'\left(L^{H^{k}_r}(r)\right)-{\bf1}_{r\le S^k_x}f_b'\left(L_k^{H^{k}_r}(r)\right)\right|.
\end{align*}
The Proposition will be proved if we show that the above right--hand side tends to zero in
$d\P\times dr$-measure, as $k\to\infty$.
Consider the first term on the right. By Corollary \ref{cor:convH},  $H^{k}_r\to H_r$ in probability, locally uniformly in~$r$, as $k\to\infty$. Moreover,  $t\to L^t(r)$ is continuous for 
$t$ outside an at most countable set, and $H$ spends zero time in that at most countable set. Hence 
we have that  $L^{H^{k}_r}(r)\to L^{H_r}(r)$ in probability, $dr$ a.e. Since $f'_b$ is continuous, the first term converges.

The second term on the r.h.s. of the previous inequality is bounded from above by
\[
\left|{\bf1}_{r\le S_x}-{\bf1}_{r\le S^k_x}\right|\, f_b'\left(L^{H^{k}_r}(r)\right)+
\left| f_b'\left(L^{H^{k}_r}(r)\right)-f_b'\left(L_k^{H^{k}_r}(r)\right)\right|.  \]
The first term in this expression converges to $0$ in $d\P\times dr$-measure thanks to Lemma \ref{le:convS}. Concerning the second term, since $f'_b$ is uniformly continuous, it suffices to show that 
\[ L^{H^{k}_r}(r)-L_k^{H^{k}_r}(r)\to 0 \quad \mbox{ in probability}.\]
Due to \eqref{defLT} and \eqref{LTk} for the local times of $H$ and $H^{k}$, this expression takes the form
\begin{align*}
&L^{H^{k}_r}(r)-L_k^{H^{k}_r}(r)=\beta(H_r-H^{k}_r)^+-\sqrt{2\beta}\int_0^r({\bf1}_{H_v>t}-{\bf1}_{H^{k}_v>t})dB_v\Big|_{t=H^{k}_r}\\
&-\int_0^r\int_{\eps_k}^\infty\!\!\!({\bf1}_{H_v>H^{k}_r}-{\bf1}_{H^{k}_v>H^{k}_r})z\tilde{N}(dv,dz)\\
&+\int_0^r\int_{\eps_k}^\infty\!\!({\bf1}_{H_v>H^{k}_r}-{\bf1}_{H^{k}_v>H^{k}_r})(z+\inf_{v\le u\le r}X_u-X_v)^+N(dv,dz)\\
&-\int_0^r\int_0^{\eps_k}\!\!{\bf1}_{H_v>t}z\tilde{N}(dv,dz)\Big|_{t=H^{k}_r}\\ &
+\int_0^r\int_0^{\eps_k}\!\!\!{\bf1}_{H_v>H^{k}_r}(z+\inf_{v\le u\le r}X_u-X_v)^+N(dv,dz).
\end{align*}
Note that we can insert the anticipative $H^{k}_r$ in the last four integrals since three of them are Stieltjes integrals, and the third is an integral with respect to a compensated Poisson point process which is independent of $H^{k}_r$.

It is plain that
\begin{align*} 
0&\le\int_0^r\int_0^{\eps_k}\!\!\!\!{\bf1}_{H_v>H^{k}_r}(z+\inf_{v\le u\le r}X_u-X_v)^+N(dv,dz)\\
&\le
\int_0^r\int_0^{\eps_k}\!\!\!\!(z+\inf_{v\le u\le r}X_u-X_v)^+N(dv,dz).
\end{align*}
Consequently from Corollary \ref{cor:estim}
\begin{align*}
\E\int_0^r\int_0^{\eps_k}\!\!\!\!{\bf1}_{H_v>H^{k}_r}(z+\inf_{v\le u\le r}X_u-X_v)^+N(dv,dz)
&\le C(r)\int_0^{\eps_k}z^2\pi(dz)\\
&\to0,\quad\text{as }k\to\infty.
\end{align*}
 In order to estimate the next to last term, we note that from Lemma \ref{le:mart}, for any given $T>0$,
\begin{align*}
\E&\left(\sup_{0\le t\le T}\left|\int_0^r\int_0^{\eps_k}\!\!\!\!\!\!{\bf1}_{H_v>t}z\tilde{N}(dv,dz)\right|^2\right)\\
&\le 2\E\left(\left|\int_0^r\int_0^{\eps_k}\!\!\!\!\!z\tilde{N}(dv,dz)\right|^2\right)
+2\E\left(\sup_{0\le t\le T}\left|\int_0^r\int_0^{\eps_k}\!\!\!\!\!\!{\bf1}_{H_v\le t}z\tilde{N}(dv,dz)\right|^2\right)
\\
&\le2r\int_0^{\eps_k}z^2\pi(dz)+8\E\left(\left|\int_0^r\int_0^{\eps_k}\!\!\!\!\!{\bf1}_{H_v\le T}z\tilde{N}(dv,dz)\right|^2\right)\\
&\le10r\int_0^{\eps_k}z^2\pi(dz)\,.
\end{align*}
Consequently
\begin{align*}
\E\left(\sup_{t\ge0}\left|\int_0^r\int_0^{\eps_k}\!\!\!\!\!\!{\bf1}_{H_v>t}z\tilde{N}(dv,dz)\right|^2\right)&=
\lim_{T\to+\infty}\E\left(\sup_{0\le t\le T}\left|\int_0^r\int_0^{\eps_k}\!\!\!\!\!\!{\bf1}_{H_v>t}z\tilde{N}(dv,dz)\right|^2\right)\\
&\le10r\int_0^{\eps_k}z^2\pi(dz)  \to 0,\quad\text{as }k\to\infty\,.
\end{align*}
We split the two previous terms into two, choosing an arbitrary $\delta>0$, which w.l.o.g. we can assume to satisfy $\delta>\eps_k$.
By the same arguments as above,
\begin{align*}\E\left|\int_0^r\int_{\eps_k}^\delta\!\!\!({\bf1}_{H_v>H^{k}_r}-{\bf1}_{H^{k}_v>H^{k}_r})(z+\inf_{v\le u\le r}X_u-X_v)^+N(dv,dz)\right| \\ \le C(r)\int_0^\delta z^2\pi(dz)
\end{align*}
and
\[\E\left(\sup_{t\ge0}\left|\int_0^r\int_{\eps_k}^\delta\!\!\!\!({\bf1}_{H_v>t}-{\bf1}_{H^{k}_v>t})z
\tilde{N}(dv,dz)\right|^2\right)\le10r\int_0^\delta z^2\pi(dz). \]
Those can be made arbitrarily small by choosing $\delta>0$ small enough.

Denoting by $N_{s,\delta}$ the number of points of $N$ in $[0,s]\times(\delta,\infty)$, we have that
\begin{align*}&\left|\int_0^r\int_\delta^\infty\!\!\!\!\!({\bf1}_{H_v>H^{k}_r}-{\bf1}_{H^{k}_v>H^{k}_r})(z+\inf_{v\le u\le r}X_u-X_v)^+N(dv,dz)\right| \\
&\qquad\le\sum_{i=1}^{N_{s,\delta}}|{\bf1}_{H_{T_i}>H^{k}_r}-{\bf1}_{H^{k}_{T_i}>H^{k}_r}|Z_i 
\end{align*}
and
\begin{align*}
&\int_0^r\int_\delta^\infty\!\!\!({\bf1}_{H_v>H^{k}_r}-{\bf1}_{H^{k}_v>H^{k}_r})  z\tilde{N}(dv,dz)\\
&=\sum_{i=1}^{N_{s,\delta}}({\bf1}_{H_{T_i}>H^{k}_r}-{\bf1}_{H^{k}_{T_i}>H^{k}_r})Z_i
-\int_{\delta}^\infty z\pi(dz)\int_0^r({\bf1}_{H_v>H^{k}_r}-{\bf1}_{H^{k}_v>H^{k}_r})dv.
\end{align*}
The fact that the finite sum converges to 0 as $k\to\infty$ follows from the fact that $H^{k}_{T_i}\to H_{T_i}$, while $H^{k}_r-H^{k}_{T_i}\to H_r-H_{T_i}\not=0$ a.s., and moreover
\[ \left|{\bf1}_{H_{T_i}>H^{k}_r}-{\bf1}_{H^{k}_{T_i}>H^{k}_r}\right|\le{\bf1}_{|H_{T_i}-H^{k}_{T_i}|>|H^{k}_r-H^{k}_{T_i}|},\]
which tends to 0 from the above claims. Moreover, by the occupation times formula,
\[ \int_0^r({\bf1}_{H_v>H^{k}_r}-{\bf1}_{H^{k}_v>H^{k}_r})dv=\int_{H^{k}_r}^\infty(L^t(r)-L^t_k(r))dt.\]
Let us use again the stopping times $\tau_K$, defined just before Lemma \ref{le:boundLT}.
Since \mbox{$\P(\tau_K\le r)$} $\to0$ as $K\to \infty$, it suffices to consider
\begin{align*} 
\E\left(\left|\int_{H^{k}_r}^\infty(L^t(r)-L^t_k(r))dt\right|;r<\tau_K\right)&\le \int_0^M\E\left(\left|L^t(r)-L^t_k(r)\right|;r<\tau_K\right) dt\\
&\quad +\int_M^\infty\E\left(\left|L^t(r)-L^t_k(r)\right|;r<\tau_K\right) dt\,.
\end{align*}
Since the integrand on the right converges to $0$ in probability for any $t$, the convergence to~$0$ of the first integral on the right follows from uniform integrability provided by Lemma~\ref{le:boundLT}, for any $M>0$. Concerning the last term, using the inequality
$\left|L^t(r)-L^t_k(r)\right|$ $\le$ \mbox{$L^t(r)+L^t_k(r)$}, we have two integrals to estimate. We estimate the first one, the estimate of the second
being very similar.
\begin{align*}
\E\int_M^\infty L^t(r){\bf1}_{r<\tau_K}dt&=\sum_{k=0}^{\infty}\E\int_{M+k}^{M+k+1}L^t(r){\bf1}_{r<\tau_K}dt\\
&\le \sum_{k=0}^{\infty}\P\left(\sup_{0\le u\le r\wedge\tau_K}H_u>M+k\right)\\
&\le \sum_{k=0}^{\infty}\P\left(\sup_{0\le u\le r\wedge\tau_K}|X_u|>\frac{\beta}{2} M\right)\\
&\le C_{2,K,r}(2/\beta)^2\sum_{k=0}^\infty (M+k)^{-2},
\end{align*}
where we have used \eqref{estimX} and Chebychev's inequality for the last line. Clearly the last right--hand side tends to $0
$ as $M\to\infty$.

Finally we consider the Brownian integral. Let us define
\[ \Phi(t,r)=\int_0^{r\wedge \tau_K}{\bf1}_{H_v>t}dB_v,\quad \Phi_k(t,r)=\int_0^{r\wedge \tau_K}{\bf1}_{H^{k}_v>t}dB_v.\]
We need to show that $\Phi(H^{k}_r,r)-\Phi_k(H^{k}_r,r)\to0$, which will follow from a variant of the last argument which we have used and the fact that for any $M>0$,
\begin{equation}\label{cvsup}
 \sup_{0\le t\le M}|\Phi(t,r)-\Phi_k(t,r)|\to0
 \end{equation}
in probability, as $k\to\infty$. It is plain that for fixed $t$, $\Phi_k(t,r)\to\Phi(t,r)$ in probability. So \eqref{cvsup} will follow if we show that for any fixed $r$, the sequence of processes $\{\Phi_k(\cdot,r)\}_{k\ge1}$ is tight in $C([0,M])$. 
It follows from the computation done in the proof of Proposition \ref{pro:contLT} and from Theorem  I.2.1 in \cite{RY} that with any
$p>2$, $\rho<\frac{1}{2}-\frac{1}{p}$, 
\[ \xi_{k,\rho}:=\sup_{0\le t\not=t'\le M}\frac{|\Phi_k(t,r)-\Phi_k(t',r)|}{|t'-t|^\rho}\]
satisfies 
\[\E[\xi_{k,\rho}^{2p}]\le C_p\E\int_0^M(L^u_k(r\wedge \tau_K))^pdu,\]
which thanks to Lemma \ref{le:boundLT} yields the desired tightness. The result follows, since 
\mbox{$\P(\tau_K< r)$} $\to 0$ as $K\to\infty$, for any $r>0$.  \epf

Let us repeat here Lemma 24 from \cite{P}.
\begin{lemma}\label{le:Gir}
Let $(\xi_k,\eta_k)$, $(\xi,\eta)$ be random pairs defined on a probability space  $(\Omega,\F,\P)$
with $\eta_k$, $\eta$ being non--negative random variables satisfying $\E[\eta_k]=\E[\eta]=1$.
Let $\tilde{\xi}_k$ stand for the r.v. $\xi_k$ defined on $(\Omega,\F,\tilde{\P}_k)$ with $d\tilde{\P_k}/d\P=\eta_k$,
$\tilde{\xi}$ for the r.v. $\xi$ defined on $(\Omega,\F,\tilde{\P})$ with $d\tilde{\P}/d\P=\eta$. If 
$(\xi_k,\eta_k)$ converges in law towards $(\xi,\eta)$ as $k\to\infty$, then $\tilde{\xi}_k$ converges in law towards
$\tilde{\xi}$, as $k\to\infty$.
\end{lemma}

Now a combination of Corollary \ref{cor:convH}, Proposition \ref{limitGirsanov} and Lemma \ref{le:Gir}
yields
\begin{proposition}\label{pro:convHakf}
As $k\to\infty$, the solution $H^{k}_s$, $s\ge 0$ of equation \eqref{eq:Hakf}
converges in probability, locally uniformly in $s$, to the solution $H_s$, $s\ge 0$, of equation \eqref{eq:Haf}.
\end{proposition}
Imitating the proof of Theorem \ref{th:RNlin}, we now deduce from Proposition~\ref{RNcpp_f} the following
\begin{theorem}\label{RNab}
 For any $a,b > 0$, under the probability measure $\P^{a,b}$, the process $\{L^t(S_x),\ 0\le t\le a,\ x>0\}$ is
a solution of the collection indexed by $x$ of SDEs \eqref{eq:Zinteractb} on the time interval $[0,a]$.
\end{theorem}

It remains to let first $b\to\infty$, then $a\to\infty$.

First of all, let us fix $x>0$.
 We would like to replace $f_b$ and $f'_b$ by $f$ and $f'$. Since $f'$ is not bounded from below, it is not clear that we can apply Girsanov's theorem, i.e. it is not clear that we have  $\E(Y^a_s)=1$ if we define
 \[Y^{a}_s=\exp\left(\frac{1}{\sqrt{2\beta}}\int_0^{s\wedge S_{x}}[f'(L^{H_r}(r))+g_a(H_r)]dB_r-\frac{1}{4\beta}\int_0^{s\wedge S_{x}} [f'(L^{H_r}(r))+g_a(H_r)]^2dr\right).\]
 We shall argue as in section 7.2.2 of \cite{P}. Let
  \[ S^n:=\inf\{s>0,\, L^{H_s}(s)>n\}\, .\]
 Since $f'(L^{H_s}(s))$ is bounded on $[0,S^n]$, we can define the probability measure $\P^a$ on $\bigvee_n\F_{S^n}$, which  is  such that for any $n\in \mathbb N$
 \[ \frac{d\P^a}{d\P}\Big|_{\F_{S^n}}=Y^a_{S^n}\, .\]
 Under $\P^a$, $H$ solves the SDE
\begin{equation}\label{eq:Ha}
\begin{split}
\beta H_s&=\int_0^s[f'(L^{H_r}(r))+g_a(H_r)]dr+X^{a}_s-\inf_{0\le r\le s}X_r\\
&\quad -\int_0^s\int_0^\infty\left(z+\inf_{r\le u\le s}X_u-X_r\right)^+N(dr,dz)\, ,
\end{split}
\end{equation}
with 
\vspace{-0.7cm}
\begin{align*}
X^a_s&=\sqrt{2\beta}B^{a}_s+\int_0^s\int_0^\infty z\tilde{N}(dr,dz),
\end{align*}
where
\vspace{-0.7cm}
\begin{align*}
B^{a}_s&=B_s-\frac{1}{\sqrt{2\beta}}\int_0^s[f'(L^{H_r}(r))+g_a(H_r)]dr
\end{align*}
is a Brownian motion under $\P^a$. It remains to verify that for each $s>0$, $\E(Y^a_s)=1$, and $\frac{d\P^a}{d\P}\Big|_{\F_s}=Y^a_s$. From Proposition 28 in \cite{P}, this will be the case, provided
\begin{lemma}
As $n\to\infty$, $\P(S^n<S_x)\to0$ and $\P^a(S^n<S_x)\to0$.
\end{lemma}
\bpf Let us establish the first statement. We choose an arbitrary $\eps>0$ and observe that for any $A>0$
\begin{align*}
\P\left(\sup_{0\le s\le S_x}L^{H_s}(s) >n\right)&\le \P\left(\sup_{0\le s\le S_x}L^{H_s}(S_x)>n\right)\\
&\le \P\left(\sup_{0\le s\le S_x}H_s>A\right)+\P\left(\sup_{0\le t\le A} L^t(S_x)>n\right)
\end{align*}
Since under $\P$, $S_x<\infty$ a.s. and $H$ is continuous, the random variable $\sup_{0\le s\le S_x}H_s$ is a.s. finite, and we can choose $A$ large enough, such that $\P\left(\sup_{0\le s\le S_x}H_s>A\right)\le \eps/2$. Next, under $\P$, the second Ray--Knight theorem 
holds true, so that
\[\P\left(\sup_{0\le t\le A} L^t(S_x)>n\right)=\P\left(\sup_{0\le t\le A} Z^x_t>n\right)\le \eps/2,\]
provided $n$ is choosen large enough.

For the proof of the second statement, we start the argument in exactly the same way. Again, thanks to the drift $g_a$,  
the random variable $\sup_{0\le s\le S_x}H_s$ is a.s. finite under $\P^a$. In order to estimate the second term, we cannot use the identification with  the solution of \eqref{eq:Zinteract} for $t>a$. However, if we  go back to the proof of Proposition \ref{RNcpp_f}, we note that for $t>a$ which had been excluded in that proof, we have
\[L^t(S_x)\le x+\int_0^{S_x}{\bf1}_{H_r\le t}f_b'(L^{H_r}(r))dr +\int_0^{S_x}{\bf1}_{H_r\le t}dB^{a}_r 
+\int_0^{S_x}{\bf1}_{H_r\le t}\int_0^\infty z\tilde{N}(dr,dz),\]
although the equality does not hold. Going through the first step of the proof of Proposition~\ref{RNcpp_f}, we deduce that $L^t(S_x)$ is a subsolution of equation \eqref{eq:Zinteract}, thus by the comparison theorem for that SDE (see \cite{DL}), $L^t(S_x)\le Z^x_t$, which finishes the proof. \epf 

It is now easy to deduce from Theorem \ref{RNab}
\begin{theorem}\label{RNa}
 For any $a>0$, under the probability measure $\P^{a}$, the process \mbox{$\{L^t(S_x),\ 0\le t\le a,\ x>0\}$} is
a solution of the collection indexed by $x$ of SDEs \eqref{eq:Zinteract} on the time interval~$[0,a]$.
\end{theorem}

It now remains to let $a\to\infty$.  First of all let us observe that the projective limit of the laws of $(H,X)$ under $\P^a$ as $a\to \infty$ renders a (unique) weak solution of \eqref{eq:H}, \eqref{eq:X}.
For {\em that} $H$ there is, however, no guarantee that $S_x<\infty$. On the other hand, the law of $\{L^t(S_x),\ 0<t\le t'\}$ depends only upon the pieces of trajectories of $H$ below $t'$, and it does not depend upon~$a$, provided $a>t'$. Therefore there exists a projective limit of those laws as well, and we have our final theorem.

 \begin{theorem}\label{RNfinal}
 There exists a random field $\{\mathcal{L}^t_x,\ x>0, t\ge0\}$ defined on a probability space $(\Omega',\F',\P')$ such that for any $a>0$, the law of $\{\mathcal{L}^t_x,\ x>0, 0\le t\le a\}$ is the same as the law of $\{L^t(S_x),\ x>0, 0<t\le a\}$ under $\P^a$. Consequently, $\{\mathcal{L}^t_x,\ x>0, t\ge0\}$ solves the collection indexed by $x>0$ of SDEs \eqref{eq:Zinteract}.
 \end{theorem}

\paragraph{Acknowledgement} The second author wishes to acknowledge useful discussions with Jean Bertoin, concerning some technical issues while writing this paper.


\begin{thebibliography}{99}
\bibitem{BP} {\sc M. Ba, E. Pardoux} Branching processes with interaction and a generalized Ray--Knight Theorem, {\it Annales de l'IHP} {\bf51}, 1290--1313, 2015.
\bibitem{BPS} {\sc M. Ba, E.  Pardoux and A.B. Sow}
 Binary trees, exploration processes, and an extended Ray--Knight Theorem, \emph{J. Appl. Probab.}  {\bf 49},  210--225, 2012.
 \bibitem{BFF} {\sc J. Berestycki, M.C. Fittipaldi, J. Fontbona} Ray--Knight representation of flows of branching processes with competition by pruning of L\'evy trees, \emph{Probab. Theory \& Relat. Fields}  {\bf172}, 725--788, 2018.
\bibitem{B} {\sc J. Bertoin} {\it L\'evy processes}, Cambridge University Press 1996.
\bibitem{DL} {\sc D. Dawson and Z. Li} Stochastic equations, flows and measure-valued processes., \emph{Ann. Probab.} {\bf 40},  813--857, 2012.
\bibitem{JFD} {\sc J.--F. Delmas}  Height process for super-critical continuous state branching process. 
\emph{Markov Proc. and Rel. Fields} {\bf 14},  309--326, 2008.
\bibitem{DP} {\sc I. Dram\'e and E. Pardoux} Approximation of a generalized CSBP with interaction, \emph{Electron. Commun. Probab.} {\bf23}, 73, 1--14, 2018.
\bibitem{DLG} {\sc T. Duquesne and J.--F. Le Gall} {\it Random trees, L\'evy processes and spatial branching processes}, Ast\'erisque {\bf281},
Soci\'et\'e Math\'ematique de France, 2002.
\bibitem{LePW} {\sc V. Le, E. Pardoux and A. Wakolbinger} Trees under attack: a Ray--Knight representation of Feller's branching diffusion with logistic growth. \emph{Probab. Theory \& Relat. Fields} {\bf155}, 583-619, 2013.
\bibitem{Li} {\sc Z. Li} {\it Measure valued branching Markov processes}, Springer Verlag 2011. 
\bibitem{P} {\sc E. Pardoux} {\it Probabilistic models of population evolution}, Springer 2016.
\bibitem{PW} {\sc E. Pardoux and A. Wakolbinger} From Brownian motion with a local time drift to Feller's branching diffusion with logistic growth,
{\it Electron. Commun. Probab.} {\bf16}, 720-731, 2011.
\bibitem{PWa} {\sc E. Pardoux and A. Wakolbinger} A path-valued Markov process indexed by the ancestral mass, {\it ALEA Lat. Am. J. Probab. Math. Stat.}{\bf 12}, 193-212, 2015.
\bibitem{PS} {\sc G. P\'olya and G. Szeg\"o} {\it Problems and Theorems in Analysis I}, Spinger, 1998.
\bibitem{Pr} {\sc P. Protter} {\it Stochastic integration and differential equations} 2d ed., Applications of Math. {\bf21}, Springer
2004.
\bibitem{RY} {\sc D. Revuz and M. Yor} {\it Continuous martingales and Brownian motion}, 3d ed., 
Grundlehren der Math. Wissenschaften {\bf293} Springer 1999.
\end{thebibliography}
\end{document}